\documentclass[12pt,a4paper]{amsart}
\usepackage{tikz}
\usepackage{pdflscape}
\usepackage{amsmath}
\usepackage{amsfonts}
\usepackage{amssymb}
\usepackage{mathtools}
\usepackage[all]{xy}
\usepackage{color}
\usetikzlibrary{arrows}
\usepackage{float}
\usepackage[shortlabels]{enumitem}
\usepackage{ifpdf}
\ifpdf
  \usepackage[colorlinks,
            linkcolor=magenta,       %%命题等式超链接颜色
            anchorcolor=magenta,     %%修改此处为你想要的颜色
            citecolor=magenta,       %%参考文献超链接颜色
            final,
            hyperindex]{hyperref}
\else
  \usepackage[colorlinks,
            linkcolor=magenta,       %%命题等式超链接颜色
            anchorcolor=magenta,     %%修改此处为你想要的颜色
            citecolor=magenta,       %%参考文献超链接颜色
            final,
            hyperindex,driverfallback=hypertex]{hyperref}
\fi
\newtheorem{theorem}{Theorem}[section]
\newtheorem{lemma}[theorem]{Lemma}
\newtheorem{coro}[theorem]{Corollary}

\newtheorem{conjecture}[theorem]{Conjecture}
\newtheorem{prop}[theorem]{Proposition}
\theoremstyle{definition}
\newtheorem{defn}[theorem]{Definition}
\newtheorem{remark}[theorem]{Remark}

\newcommand{\delete}[1]{}

\newcommand {\Des}{{\rm{Des}}}

\newcommand {\Inv}{{\rm{Inv}}}
\newcommand {\Nsp}{{\rm{Nsp}}}

\newcommand {\inv}{{\rm{inv}}}
\newcommand {\nsp}{{\rm{nsp}}}

\newcommand {\oinv}{{\rm{oinv}}}

\newcommand {\onsp}{{\rm{onsp}}}
\newcommand {\sgn}{{\rm{sgn}}}

\newcommand\barr[1]{\overline{#1}}
\begin{document}
\title[Sign-twisted generating functions of the odd length]{Sign-twisted generating functions of the odd length for Weyl groups of type $D$}

\author{Haihang Gu}
\address{School of Mathematics and Statistics, Southwest University, Chongqing 400715, China;}
\address{School of Mathematical Sciences, East China Normal University, Shanghai 200241, China}
\email{52275500014@stu.ecnu.edu.cn}

\author{Houyi Yu$^{\ast}$}\thanks{*Corresponding author}
\address{School of Mathematics and Statistics, Southwest University, Chongqing 400715, China}
\email{yuhouyi@swu.edu.cn}

%\date{\today}

\begin{abstract}
The odd length in Weyl groups is a new statistic analogous to the classical Coxeter length, and features combinatorial and parity conditions.
We establish explicit closed product formulas for the sign-twisted generating functions of the odd length for parabolic quotients of
Weyl groups of type $D$.
As a consequence, we verify three conjectures of Brenti and Carnevale on evaluating closed forms for these generating functions.
We then give an equivalent condition for the sign-twisted generating functions to be expressible as products of cyclotomic polynomials, settling a conjecture of Stembridge.
%We show that there is at most one non-cyclotomic factor of the form $x^n + 2x^m + 1$ for some positive integers $m$ and $n$ in such a sign-twisted generating function.
\end{abstract}


\maketitle

%\vspace{-0.8cm}

\tableofcontents

\vspace{-0.8cm}

\section{Introduction}

\subsection{Odd length}
The study of statistics on finite Coxeter groups has been a significant problem in algebra and combinatorics for the last three decades.
This is due partly to the fact that a better understanding of various numerical statistics and finding their (signed) generating functions
is one effective approach to look at these groups. For example, the generating function of the Coxeter length for a finite Coxeter group encodes all degrees of the group \cite{Hum90}.
There has been tremendous interest in developing signed enumeration of finite Coxeter groups by various statistics (see, for instance, \cite{AGR05, Bia06, Cas12, DF92, Mon15, Rei95,Siv11,Wac92} and the references therein).

Recently, motivated by questions in enumerative geometry, a new statistic called the odd length was introduced on Weyl groups.
This statistic is analogous to the classical Coxeter length, and combines combinatorial and parity conditions.
More precisely, given a Weyl group $W$ with a crystallographic root system,
for any $w\in W$, the \emph{odd length} of $w$,
denoted $L(w)$, is the number of positive roots with odd height sent to negative roots by $w$.

The concept of the odd length was first introduced in Weyl groups of type $A$ by Klopsch and Voll in \cite{KV09}, where it was used to deal with the enumeration of non-degenerate flags in  finite vector spaces equipped with various classes of bilinear forms.
Then the odd length in Weyl groups of type $B$, as a natural analogue of that in type $A$, was defined by Stasinski and Voll \cite{SV14} in relation to the local factors of the representation zeta functions of finitely generated torsion-free nilpotent groups.
An analogous statistic has later been defined and studied in Weyl groups of type $D$ in \cite{BC17D}.
By abstracting the common features of the odd lengths in the classical Weyl groups,
Brenti and Carnevale in \cite{BC19} defined, for any crystallographic root system, an analogous odd length statistic in the corresponding Weyl group.
%This new statistic depends on the choice of the simple system in the root system, whereas, its sign-twisted generating function over the whole Weyl group only depends on the root system.
By choosing a natural simple system, it agrees exactly, for the classical Weyl groups, with the definition of the odd length that investigated in \cite{BC17D,KV09,SV13,SV14}.
Developing various combinatorial properties of the odd length is an active area of research in recent years
(see, for example, \cite{BC21, BCT23, BC20, FG23,Lan18,Ste19}).

The main purpose of this paper is to establish closed product formulas for the sign-twisted generating functions of the odd length for parabolic quotients of Weyl groups of type $D$, thereby verifying several conjectures on these generating functions.

\subsection{Sign-twisted generating functions}
Given a Weyl group $(W,S)$,
for each $I\subseteq S$, we write $W_I$ for the standard \emph{parabolic subgroup} of $W$ generated by $I$, and $W^I$ for the unique \emph{left cosets representatives of minimum length} of $W_I$ in $W$.
The \emph{sign-twisted  generating function} of the odd length over $W^I$ is defined by
\begin{align*}%\label{eq:SGFOL}
W^I(x):=\sum_{w\in W^I}(-1)^{\ell(w)}x^{L(w)},
\end{align*}
where $\ell(w)$ is the usual Coxeter length of $w$.
A remarkable property of the odd length is that its sign-twisted
generating function over the corresponding Weyl group always has an elegant factorization \cite{BC19,Ste19}.

The Weyl groups that we focus on primarily in this paper are of types $A_{n-1}$, $B_n$ and $D_n$, which can be expressed combinatorially as the symmetric group $\mathcal{S}_n$, the hyperoctahedral group $\mathcal{B}_n$ and the even hyperoctahedral group $\mathcal{D}_n$, respectively.
For a suitable set of generators, the simple system of type $A_{n-1}$ can be identified with the set $\{1,\ldots,n-1\}$, while the simple systems of types $B_n$ and $D_n$ can be identified with the set $\{0,1,\ldots,n-1\}$.

Explicit closed product formulas of $\mathcal{S}_n^I(x)$ and $\mathcal{B}_n^I(x)$ were conjectured in \cite[Conjecture C]{KV09} and \cite[Conjecture 1.6]{SV14}, respectively.
These conjectures were verified by Stasinski and Voll \cite{SV13} in type $B$ for several special cases, and fully proved independently by Brenti and Carnevale \cite{BC17T} in types $A$ and $B$, and by Landesman \cite{Lan18} in type $B$ in a different way.
For descent classes which are singletons the formula in type $B$ is equivalent to the Poincar\'e polynomial of the varieties of symmetric matrices of fixed rank over a finite field; see, for instance, \cite[Section 1]{SV13}, or \cite[Lemma 10.3.1]{Igu00} and compare with \cite[Lemma 3.1(3.4)]{SV14}.

To state our main results we introduce some notation.
Let $\mathbb{P}$ denote the set of positive integers and $\mathbb{N}:=\mathbb{P}~\cup \{0\}$.
Given any integers $a$ and $b$, let $[a,b]:=\{a,a+1,\ldots,b\}$ if $a\leq b$, and $[a,b]:=\emptyset$ if $a> b$.
In this paper, $n$ denotes a nonnegative integer.
For notational simplicity, we usually write $\overline{n}$ for the additive inverse $-n$ of $n$.
Set $[n]:=[1,n]$ and $[n]_{\pm}:=[n]\cup[\overline n,\overline1]$.
For a real number $x$ we denote by $\lfloor x\rfloor$ the greatest integer less than or equal to $x$, and similarly by $\lceil x\rceil$ the smallest integer greater than or equal to $x$.
Given a set $I$ we denote by $|I|$ its cardinality.
%If $I=\{a_1,a_2,\ldots,a_n\}$ with $a_1<a_2<\cdots<a_n$, then we write $I=\{a_1,a_2,\ldots,a_n\}_{<}$ for simplicity.

For each subset $I$ of $[0,n-1]$, there exist uniquely $a_0<a_1<\cdots<a_s$ and $b_1<b_2<\cdots<b_s$ in $[0,n-1]$
such that $I=I_0\cup I_1\cup\cdots\cup I_s$, where $I_0:=[0, a_0-1]$ and $I_i:=[a_i, b_i]$ for all $i\in[s]$.
Notice that here we have $b_{i-1}+1<a_i\leq b_i$ for all $i\in[s]$, where $b_0:=a_0$.
These intervals $I_i$, $i\in[0,s]$, are called the \emph{connected components} of $I$.
Then $I_0$ is the unique connected component containing $0$. In particular, $I_0=\emptyset$ if $|I_0|=0$.

For nonnegative integers $n_1,n_2,\ldots,n_k$ summing to $n$, we define the \emph{$q-$multinomial coefficient}
$$
\left[\begin{matrix}n\\ n_1,n_2,\ldots,n_k\end{matrix}\right]_q:=\frac{[n]_q!}{[n_1]_q![n_2]_q!\cdots[n_k]_q!},
$$
where
$$
[n]_q:=\frac{1-q^n}{1-q}, \quad [n]_q!:=\prod_{i=1}^{n}[i]_q, \quad \text{and}\quad [0]_q!:=1.
$$
Given a subset $I\subseteq [0,n-1]$ with connected components $I_0, I_1,\ldots, I_s$, we write
\begin{align}\label{eq:m_I}
m_I:=\sum_{k=0}^s\left\lfloor\frac{|I_k|+1}{2}\right\rfloor
\quad \text{and} \quad
C_I:=\left[\begin{matrix}m_I\\ \left\lfloor\frac{|I_0|+1}{2}\right\rfloor,\left\lfloor\frac{|I_1|+1}{2}\right\rfloor,\ldots,\left\lfloor\frac{|I_s|+1}{2}\right\rfloor\end{matrix}\right]_{x^2}.
\end{align}

Following \cite{BC17D}, the formulas of $\mathcal{S}_n^I(x)$ and $\mathcal{B}_n^I(x)$
can be stated as follows. %Here we refer the reader to Section \ref{sec:Prel} for more detailed notions.
\begin{theorem}[\cite{BC17D}, Theorem 2.8]\label{theo:SnI(x)}
Let $I\subseteq[n-1]$. Then
\begin{align*}
\mathcal{S}_n^I(x)=C_I\prod_{j=2m_I+2}^n\left(1+(-1)^{j-1}x^{\left\lfloor\frac{j}{2}\right\rfloor}\right).
\end{align*}
\end{theorem}
Notice that the $n$ in Theorem \ref{theo:SnI(x)} may be less than $2m_I+2$. For example, if $n=8$ and $I=[1,3]\cup[5,7]$, then $m_I=4$ and hence $n=2m_I$.
We will assume throughout that an empty product is equal to $1$.

\begin{theorem}[\cite{BC17D}, Theorem 2.10]\label{theo:BnI(x)1}
Let $I\subseteq[0,n-1]$ and let $I_0,I_1,\ldots,I_s$ be the connected components of $I$. Then
\begin{align*}
\mathcal{B}_n^I(x)=\left[\begin{matrix}m\\ \left\lfloor\frac{|I_1|+1}{2}\right\rfloor,\ldots,\left\lfloor\frac{|I_s|+1}{2}\right\rfloor\end{matrix}\right]_{x^2}\frac{\prod_{j=|I_0|+1}^n(1-x^j)}{\prod_{i=1}^m(1-x^{2i})},
\end{align*}
where $m:=\sum_{k=1}^s\left\lfloor\frac{|I_k|+1}{2}\right\rfloor$.
\end{theorem}

In \cite{BC17D}, Brenti and Carnevale %studied the the odd length and its sign-twisted generating functions on the Weyl groups of type $D$. They
showed that many properties satisfied by $\mathcal{S}_n^I(x)$ and $\mathcal{B}_n^I(x)$ still hold for $\mathcal{D}_n^I(x)$.
In particular, closed product formulas of $\mathcal{D}_n^I(x)$ were obtained when $|I|\leq 1$, $I=\{0,1\}$ and $I=\{0,2\}$.
Based on these formulas have been established, Brenti and Carnevale proposed the following conjectures, revealing explicit evaluations of $\mathcal{D}_n^I(x)$.

\begin{conjecture}[\cite{BC17D}, Conjectures 6.1 and 6.2]\label{conj:DnI0i}
If $n\geq 5$ and $3\leq i\leq n-1$, Then
\begin{align*}
\mathcal{D}_n^{\{0,i\}}(x)=\prod_{j=4}^n\left(1+(-1)^{j-1}x^{\left\lfloor\frac{j}{2}\right\rfloor}\right)^2\quad\text{and}\quad
\mathcal{D}_n^{\{0,1,i\}}(x)=(1-x^4)\prod_{j=5}^n\left(1+(-1)^{j-1}x^{\left\lfloor\frac{j}{2}\right\rfloor}\right)^2.
\end{align*}
\end{conjecture}
\begin{conjecture}[\cite{BC17D}, Conjecture 6.5]\label{conj:DnI(x)}
Let $n\geq 3$, and let $I\subseteq [0,n-1]$ with connected components $I_0, I_1,\ldots, I_s$. Then there exists a polynomial $M_I(x)\in \mathbb{Z}[x]$ such that
\begin{align*}
\mathcal{D}_n^{I}(x)=M_I(x)\prod_{j=2m_I+2}^n \left(1+(-1)^{j-1}x^{\left\lfloor\frac{j}{2}\right\rfloor}\right)^2.
\end{align*}
Furthermore, $M_I(x)$ only depends on $(|I_0|, |I_1|,\ldots, |I_s|)$ and is a symmetric function of $|I_1|,\ldots, |I_s|$.
\end{conjecture}
%Conjecture \ref{conj:DnI0i} was formulated as \cite[Conjectures 6.3 and 6.4]{BC17D} in an equivalent way in terms of the sign-twisted generating functions in type $A$.
Conjecture \ref{conj:DnI(x)} has been verified for the whole group $\mathcal{D}_n$ and its maximal, and some other, quotients.
We will give an explicit closed product expression for $\mathcal{D}_n^{I}(x)$ (see Theorem \ref{theo:finaltheorem}), resolving Conjectures \ref{conj:DnI0i} and \ref{conj:DnI(x)}.
It turns out that the polynomial $M_I(x)$ always exists, but its form in general depends on $(I_0, I_1,|I_2|,\ldots, |I_s|)$, and $M_I(x)$ is symmetric only with respect to $|I_2|,\ldots, |I_s|$.

A crucial property of $\mathcal{S}_n^{I}(x)$ and $\mathcal{B}_n^{I}(x)$ is that they are products of cyclotomic polynomials.
Stembridge \cite{Ste19} conjectured this is also true nearly for all $\mathcal{D}_n^{I}(x)$.
\begin{conjecture}[\cite{Ste19}, Conjecture 7.4(a)]\label{conj:Ste197.4a}
Let $W=W(\Phi)$ be a Weyl group, where $\Phi$ is a crystallographic irreducible root system of type $D$ with a simple system $\Delta$. Then for any $I\subseteq\Delta$, $W^I(x)$ is a product of  cyclotomic polynomials unless
$\Phi\cong D_{2m}$ for some positive integer $m$ and $K\subseteq I\subseteq \Delta$,
where $\Phi_K\cong A_1^{\oplus m+1}$.
\end{conjecture}
We give a necessary and sufficient condition for $\mathcal{D}_n^I(x)$ that can not be expressed as a product of cyclotomic polynomials, thereby confirming Conjecture \ref{conj:Ste197.4a} (see Corollary \ref{coro:Stembconj}).
In particular, we show that the $\mathcal{D}_n^I(x)$ that is not a product of cyclotomic polynomials contains exactly one non-cyclotomic factor, which always has the form $x^n + 2x^m + 1$ for some positive integers $m$ and $n$.

\subsection{Outline}

We explain below our approach to proving the main result Theorem \ref{theo:finaltheorem} and at the same time give the organization of the paper.

After summarizing some necessary concepts and basic facts on Coxeter systems, we give in Section \ref{sec:Prel} the definition of the odd length in Weyl groups and its combinatorial description in terms of permutations for the root systems of types $A$  and $D$.

Section \ref{sec:I0neqemptyset} is devoted to the study of some general properties of the sign-twisted generating functions $\mathcal{D}_n^I(x)$.
The concept of a \emph{compressed subset} of $[0,n-1]$ (see Definition \ref{defn:compressedset})
is introduced, generalizing the one in \cite{BC17T} defined only for subsets of $[n-1]$.
We show that for each subset $I$ of $[0,n-1]$ that contains $0$ and $1$, the computation of $\mathcal{D}_n^I(x)$ can be reduced to the case where $I$ is a compressed subset of the form $I=[0,a_0-1]\cup[a_0+1,a_1-1]\cup\cdots\cup[a_{s-1}+1,n-1]$ (see Propositions \ref{prop:a0=odd} and \ref{prop:a0=even}).
A key role for these proofs is the use of sign-reversing involutions to eliminate cancellations in an alternating sum.

In Section \ref{sec:Ceos}, we show that in order to determine $\mathcal{D}_n^I(x)$ when $I=[0,a_0-1]\cup[a_0+1,a_1-1]\cup\cdots\cup[a_{s-1}+1,n-1]$ is compressed and $a_0\geq2$,
it suffices to assume that $I=[0,a_0-1]\cup[a_0+1,n-1]$ (see Proposition \ref{lem:DnI=DnJxSn-aKeq}).
%This is achieved in two steps. First, following the terminologies of \cite{KV09,SV13}, we define \emph{chessboard elements} in the Weyl group $\mathcal{D}_n$
%(see Definition \ref{defn:chessboardelem}), and show that $\mathcal{D}_n^I(x)$ is supported on the set of chessboard elements. Second, we introduce the notion of \emph{odd sandwiches} for even signed permutations (see Definition \ref{def:oddsanwichin[k+1,n]}), extending the notion defined for signed permutations in \cite{SV13}. Then we show that $\mathcal{D}_n^I(x)$ is indeed supported on a set of special chessboard elements defined in terms of odd sandwiches.
%A crucial property of this relatively small supporting set is that the function $L$ on it is additive with respect to a suitable parabolic factorization.
A crucial point in the proof is that we must find an appropriate supporting set for $\mathcal{D}_n^{I}(x)$ on which the odd length function $L$ is additive with respect to a suitable parabolic factorization.

Having the necessary ingredients at hand, we establish an explicit closed product formula for $\mathcal{D}_n^I(x)$ in Section \ref{sec:formudnixfull} (see Theorem \ref{theo:finaltheorem}), and verify Conjectures \ref{conj:DnI0i} and \ref{conj:DnI(x)} as corollaries.

Finally, in Section \ref{sec:cyclotomicpoly}, we give combinatorial descriptions of sets of indices $I$ such that $\mathcal{D}_n^I(x)$  can not be expressed as a product of cyclotomic polynomials, proving Conjecture \ref{conj:Ste197.4a}.

\section{Preliminaries}\label{sec:Prel}

We begin by reviewing some notation, definitions and facts (mostly following  \cite{BB05, BC17D, BC19, Hum90}) that will be used throughout the paper.

%Recall that a Coxeter system is a $(W,S)$ tuple. $S$ is the set of generators of the group $W$ subject to relations of the form
%\begin{align*}
%\end{align*}
%where $s,s'\in S$, $e$ is the identity element of $W$, and $m(s,s')=1$ if and only if $s=s'$. If there is no relation between $s$ and $s'$, we denote $m(s,s')=\infty$. $W$ is called a Coxeter group and $S$ is the set of Coxeter generators.
%We follow \cite{BB05} for the notation and terminology about Coxeter group.
Let $\Phi$ be a finite crystallographic root system with a set of simple roots $\Delta$. Let $\Phi^+$ denote the corresponding set of roots that are nonnegative linear combinations of simple roots, and $\Phi^{-}:=-\Phi^{+}$, so that $\Phi=\Phi^+\cup \Phi^-$.
The group $W$ generated by the set of reflections $\{s_{\alpha}\,|\, \alpha\in \Phi\}$ is the \emph{Weyl group} of $\Phi$.
Moreover, $(W,S)$ is a Coxeter system, where $S:=\{s_{\alpha}\,|\, \alpha\in \Delta\}$.
%We let $\ell_{\Delta}$ be the \emph{Coxeter length function} on $W$ with respect to $S$.
For an element $w=s_{\alpha_1}s_{\alpha_2}\cdots s_{\alpha_r}\in W$ and a root $\alpha\in \Phi$, we let $w(\alpha)$ denote the \emph{action} of $w$ on $\alpha$ as a composition of the reflections $s_{\alpha_1},s_{\alpha_2},\ldots, s_{\alpha_r}$.
The \emph{Coxeter length} of $w$, denoted $\ell(w)$, is the number of positive roots sent to negative
roots by $w$, that is,
$
\ell(w):=|\{\alpha\in \Phi ^+\,|\, w(\alpha)\in \Phi^{-}\}|.
$
For each $I\subseteq S$, let $W_I$ denote the \emph{parabolic subgroup} generated by $I$, and define the \emph{parabolic quotient} of $W_I$ in $W$ to be
$$
W^I:=\{w\in W\,|\,\ell(ws)>\ell(w)\ \text{for all}\ s\in I\}.
$$
It is well known that any $w\in W$ has a unique factorization $w=w^I w_I$ such that $w^I\in W^I$ and $w_I\in W_I$. Moreover, for this factorization,
we have $\ell(w)=\ell(w^I)+\ell(w_I)$.

Given a root $\alpha$ in $\Phi$, write uniquely $\alpha=\sum_{\beta\in \Delta}c_\beta\beta$,  we call the sum $\sum_{\beta\in\Delta }c_\beta$ the \emph{height} of $\alpha$ (relative to $\Delta$), abbreviated ${\rm ht}(\alpha)$.
For any $w\in W$, the \emph{odd length} of $w$, denoted $L(w)$, is the number of positive roots with odd height mapped to negative roots by $w$. In other words, we have
$$
L(w)=|\{\alpha\in \Phi ^+\,|\, {\rm ht}(\alpha)\equiv1(\bmod\, 2)\ \text{and}\ w(\alpha)\in \Phi^{-}\}|.
$$

Throughout this paper, we only deal with results for crystallographic irreducible root systems of types $A$ and $D$. Our conventions for realizations of these root systems are below:
\begin{itemize}
\item Type $A_{n-1}$: $\Phi=\{\pm(e_i-e_j)\,|\,1\leq i<j\leq n\}$, $\Delta=\{e_{i+1}-e_i\,|\, i\in[n-1]\}$;
\vspace{2mm}
%\item Type $B_{n}$: $\Phi=\{\pm(e_i\pm e_j)\,|\,1\leq i<j\leq n\}\cup \{\pm e_i\,|\,i\in [n]\}$ with simple system $\Delta=\{e_1\}\cup\{e_{i+1}-e_i\,|\, i\in[n-1]\}$.
%\vspace{2mm}
\item Type $D_{n}$: $\Phi=\{\pm(e_i\pm e_j)\,|\,1\leq i<j\leq n\}$, $\Delta=\{e_1+e_2\}\cup\{e_{i+1}-e_i\,|\, i\in[n-1]\}$.
\end{itemize}

As discussed in \cite[Chapters 1 and 8]{BB05}, with the above choices of simple systems, the Weyl groups of types $A$ and $D$ have
combinatorial descriptions as symmetric groups and even hyperoctahedral groups, respectively.
Recall that a \emph{signed permutation} $\sigma$ of degree $n$ is a bijection of the set $[n]_{\pm}$ onto itself such that $\sigma(\overline i)=\overline{\sigma(i)}$ for all $i\in[n]$.
The set of all signed permutations of degree $n$ naturally form a group under composition, called the \emph{hyperoctahedral group} and denoted by $\mathcal{B}_n$, which is a Coxeter group of type $B_n$.
For a signed permutation $\sigma$ of degree $n$, we use both the one-line notation $\sigma=\sigma(1)\sigma(2)\cdots \sigma(n)$ and the disjoint cycle notation.
The set of all signed permutations of degree $n$ having an even number of negative entries in their one-line notation naturally forms the \emph{even hyperoctahedral group} which denoted in the literature by $\mathcal{D}_n$. Indeed, the group $\mathcal{D}_n$ is a Coxeter group of type $D_n$ with the generating set $S^D:=\{s_0^D, s_1^D,\ldots, s_{n-1}^D\}$, where $s_0^D:=(1,\overline 2)(2,\overline 1)$ and $s_i^D:=(i,i+1)(\barr i,\overline{i+1})$ for all $i\in[n-1]$.

The \emph{symmetric group} $\mathcal{S}_n$ is the subgroup of $\mathcal{D}_n$ consisting of all permutations of the set $[n]$,
and it is a Coxeter group of type $A_{n-1}$ with respect to the
set of Coxeter generators $S^A:=\{s_1,\ldots,s_{n-1}\}$, where $s_i:=(i,i+1)$ for all $i\in[n-1]$.

In the above combinatorial interpretations of the classical Weyl groups, we identify the generating set $S^D$ of $\mathcal{D}_n$ with the set $[0,n-1]$ via $s_i^D\leftrightarrow i$,
and similarly  identify $S^A$ with $[n-1]$.
%If $\sigma\in \mathcal{D}_n$, the descent set of $\sigma$ of type $D$  is $\{i\in[0,n-1]\,|\,\sigma(i)>\sigma(i+1)\}$, where $\sigma(0)=\overline{\sigma(2)}$.
For any $I\subseteq [0,n-1]$, we have
\begin{align*}
\mathcal{D}_n^I=\{\sigma\in \mathcal{D}_n\,|\, \sigma(i)<\sigma(i+1)\ \text{for all}\ i\in I,\ \text{where} \ \sigma(0):=\overline{\sigma(2)}\}.
\end{align*}
The \emph{sign-twisted  generating function} of the odd length over $\mathcal{D}_n^I$ is
\begin{align*}%\label{eq:SGFOL}
\mathcal{D}_n^I(x):=\sum_{w\in \mathcal{D}_n^I}(-1)^{\ell(w)}x^{L(w)}.
\end{align*}

The Coxeter length and odd length also have combinatorial interpretations in terms of statistics on the one-line notation.
Let $\sigma$ be a signed permutation of degree $n$.
Define %$\inv(\sigma):=|\Inv(\sigma)|$, $\nega(\sigma):=|\Neg(\sigma)|$, $\nsp(\sigma):=|\Nsp(\sigma)|$, where
\begin{align*}
\inv(\sigma):=&\ |\{(i,j)\in [n]^2\,|\,i<j, \sigma(i)>\sigma(j)\}|,\\
%\nega(\sigma):=&\ |\{i\in [n]\,|\,\sigma(i)<0\}|,\\
\nsp(\sigma):=&\ |\{(i,j)\in [n]^2\,|\,i<j, \sigma(i)+\sigma(j)<0\}|.
\end{align*}
Similarly, define their ``odd'' analogues %$\oinv(\sigma):=|\OInv(\sigma)|$, $\oneg(\sigma):=|\ONeg(\sigma)|$ and $\onsp(\sigma):=|\ONsp(\sigma)|$, where
\begin{align*}
\oinv(\sigma):=&\ |\{(i,j)\in [n]^2\,|\,i<j, \sigma(i)>\sigma(j), i\not\equiv j(\bmod\,2)\}|,\\
%\oneg(\sigma):=&\ |\{i\in [n]\,|\,\sigma(i)<0, i\equiv1(\bmod\,2)\}|,\\
\onsp(\sigma):=&\ |\{(i,j)\in [n]^2\,|\,i<j, \sigma(i)+\sigma(j)<0), i\not\equiv j(\bmod\,2)\}|.
\end{align*}
By \cite[Propositions 2.3 and 3.4]{BC19}, for any $\sigma\in W(\Phi)$,
\begin{itemize}
\item if $\Phi$ is of type $A$, then $\ell(\sigma)=\inv(\sigma)$ and $L(\sigma)=\oinv(\sigma)$;
%\item if $\Phi$ is of type $B$, then $\ell(\sigma)=\inv(\sigma)+\nega(\sigma)+\nsp(\sigma)$ and $L(\sigma)=\oinv(\sigma)+\oneg(\sigma)+\onsp(\sigma)$;
\item if $\Phi$ is of type $D$, then $\ell(\sigma)=\inv(\sigma)+\nsp(\sigma)$ and $L(\sigma)=\oinv(\sigma)+\onsp(\sigma)$.
\end{itemize}
For example, if $\sigma=3\barr251\barr4$, then for the lengths in type $D$, we have $\ell(\sigma)=11$ and $L(\sigma)=7$.

Since if $\sigma\in \mathcal{S}_n$, then $\sigma(i)>0$ for all $i\in[n]$, we see that the restriction of $L$ in type  $D$ to $\mathcal{S}_n$ agree with the function $L$ in type $A$.
\delete{Because we concern mainly with distributions in type $D$, we simply denote by $\ell$ and $L$ the length and odd length, respectively, on $\mathcal{D}_n$,
and similarly for $\Des$, if there exists no confusion.}
For future applications, we give the following result about $\ell$ and $L$ on $\mathcal{D}_n$.

\begin{lemma}\label{lem:ellLsigma'a=n}
Let $a\in[n]$ and $\sigma\in \mathcal{D}_n$ with $\sigma(a)>0$. If $\sigma'=(1,\overline{1})(\sigma(a),\overline{\sigma(a)})\sigma$, then
\begin{align}
	\ell(\sigma')&=\ell(\sigma)+2|\{i\in[a-1]\,|\,|\sigma(i)|<\sigma(a)\}|,\label{eq:ell(sigma)'} \\
L(\sigma')&= L(\sigma)+2|\{i\in[a-1]\,|\,|\sigma(i)|<\sigma(a),\ i\not\equiv a(\bmod\,2)\}|.\label{eq:L(sigma)'}
\end{align}
In particular, if $\sigma(a)=n$, then
$
	\ell(\sigma')=\ell(\sigma)+2a-2$ and $L(\sigma')= L(\sigma)+2\left\lfloor\frac{a}{2}\right\rfloor.
$
\end{lemma}
\begin{proof}
Since $\nsp(\sigma')=\nsp(\sigma)+|\{i\in[a-1]\cup[a+1,n]\,|\,|\sigma(i)|<\sigma(a)\}|$ and
\begin{align*}
	\inv(\sigma')=\inv(\sigma)+|\{i\in[a-1]\,|\,|\sigma(i)|<\sigma(a)\}|-|\{i\in[a+1,n]\,|\,|\sigma(i)|<\sigma(a)\}|,
\end{align*}
it follows that \eqref{eq:ell(sigma)'} and \eqref{eq:L(sigma)'} hold.
If $\sigma(a)=n$, then $|\{i\in[a-1]\,|\,|\sigma(i)|<\sigma(a)\}|=a-1$ and
\begin{align*}
|\{i\in[a-1]\,|\,|\sigma(i)|<\sigma(a),\ i\not\equiv a(\bmod\,2)\}|=\left\lfloor\frac{a}{2}\right\rfloor,
\end{align*}
so the proof follows.
\end{proof}

\section{Shifting and compressing}\label{sec:I0neqemptyset}

In this section, we recall and develop some important properties of $\mathcal{D}_n^I(x)$.
In particular, we prove that some connected components of $I$ can be compressed and shifted to the left or to the right, without changing the value of $\mathcal{D}_n^I(x)$.
The main results of this section are Propositions \ref{prop:a0=odd} and \ref{prop:a0=even}, which reduce the problem of determining $\mathcal{D}_n^I(x)$ for general $I$ to the case where $I$ is a special compressed subset of $[0,n-1]$.

\delete{
The main approach we use to deal with $\mathcal{D}_n^I(x)$ is the sign-reversing involutions, which is
one of the standard combinatorial techniques for eliminating cancellations.
Recall that  for a set $A$ equipped with a sign function ${\rm sgn}:A\rightarrow\{1,-1\}$, an involution $\varphi:A\rightarrow A$ is \emph{sign reversing} if
${\rm sgn}\,\varphi(a)=-{\rm sgn}\, a $ when $\varphi(a)\neq a$. It follows that
$$
\sum_{a\in A}{\rm sgn}a=\sum_{a\in F}{\rm sgn}a,
$$
where $F$ is the set of fixed points of $\varphi$. In particular, if $\varphi$ has no fixed points, then $\sum\limits_{a\in A}{\rm sgn}a=0$.
%Such a method of proof is called a \emph{sign-reversing involution argument}.
}
\subsection{Shifting}\label{subsec:shifting}

The following two conclusions imply that a connected component of a subset $I$ can be shifted to the left or to the right in some situations, leaving the corresponding sign-twisted generating function unchanged.
%up to multiplication by a nonzero polynomial.
\begin{lemma}[\cite{BC17D}, Proposition 3.4]\label{prop:D_n^Icup0=Icup1}
Let $I\subseteq[2,n-1]$, where $n\geq 2$. Then $\mathcal{D}_n^{I\cup \{0\}}(x)=\mathcal{D}_n^{I\cup \{1\}}(x)$.
\end{lemma}

\begin{lemma}\label{prop:shift[0,n-1]}
Let $I\subseteq[0,n-1]$, and let $\widetilde I:=(I\setminus\{i\})\cup\{i+2k+1\}$, where $k\in\mathbb{N}$ and $i\in\mathbb{P}$ with $i\geq3$.
If $i+2k+2\in [n]\backslash I$ and $[i,i+2k]$ is a connected component of $I$, then
\begin{enumerate}
\item\label{eq:prop3.7dnx0} $\mathcal{D}_{n}^I(x)=\mathcal{D}_{n}^{I\cup \widetilde{I}}(x)=\mathcal{D}_{n}^{\widetilde{I}}(x)$,\vspace{2mm}
\item\label{eq:prop3.7dnx0a=n} for any $b\in [n]\backslash[i,i+2k+2]$ and $n'\in \{n,\overline{n}\}$, we have
\begin{align*}
\sum_{\substack{\sigma\in \mathcal{D}_n^I,\, \sigma(b)=n'}}(-1)^{\ell(\sigma)}x^{L(\sigma)}
=\sum_{\substack{\sigma\in \mathcal{D}_n^{I\cup \widetilde I},\, \sigma(b)=n'}}(-1)^{\ell(\sigma)}x^{L(\sigma)}
=\sum_{\substack{\sigma\in \mathcal{D}_n^{\widetilde I},\, \sigma(b)=n'}}(-1)^{\ell(\sigma)}x^{L(\sigma)}.
\end{align*}
\end{enumerate}
\end{lemma}
\begin{proof}
\eqref{eq:prop3.7dnx0} This is \cite[Proposition 3.7]{BC17D}.

\eqref{eq:prop3.7dnx0a=n}
Since $i\geq3$ and $[i,i+2k]$ is a connected component of $I$, we have $\sigma(i)<\sigma(i+1)<\cdots<\sigma(i+2k+1)$ for any $\sigma\in \mathcal{D}_n^I$.
Then, by analyzing the value of $\sigma(i+2k+2)$, we have
\begin{align*}
\sum_{\substack{\sigma\in \mathcal{D}_n^I,\, \sigma(b)=n'}}(-1)^{\ell(\sigma)}x^{L(\sigma)}
=&\, \sum_{\substack{\sigma\in \mathcal{D}_n^I,\, \sigma(b)=n'\\ \sigma(i+2k+2)<\sigma(i)}}(-1)^{\ell(\sigma)}x^{L(\sigma)}
+\sum_{\substack{\sigma\in \mathcal{D}_n^I,\, \sigma(b)=n'\\ \sigma(i+2k+1)<\sigma(i+2k+2)}}(-1)^{\ell(\sigma)}x^{L(\sigma)}\\
&\hspace{8mm}+\sum_{j=1}^{2k+1}\sum_{\sigma\in A_j}(-1)^{\ell(\sigma)}x^{L(\sigma)},
\end{align*}
where
$A_j=\{\sigma\in \mathcal{D}_n^I\mid \sigma(b)=n', \sigma(i+j-1)<\sigma(i+2k+2)<\sigma(i+j)\}$ for $j\in[2k+1]$.
Notice that
\begin{align*}
\sum_{\substack{\sigma\in \mathcal{D}_n^I,\, \sigma(b)=n'\\ \sigma(i+2k+2)<\sigma(i)}}(-1)^{\ell(\sigma)}x^{L(\sigma)}
+\sum_{\sigma\in A_1}(-1)^{\ell(\sigma)}x^{L(\sigma)}
=\sum_{\sigma\in A_1}\left((-1)^{\ell(\sigma)}x^{L(\sigma)}+(-1)^{\ell(\sigma')}x^{L(\sigma')}\right),
\end{align*}
where $\sigma'=(\sigma(i),\sigma(i+2k+2))(\overline{\sigma(i)},\overline{\sigma(i+2k+2)})\sigma$.
A straightforward computation shows that $\ell(\sigma')=\ell(\sigma)+1$ and $L(\sigma')=L(\sigma)$,
so the above sum is zero. Similarly, for each $j\in[2, 2k+1]$ with $j\equiv 0(\bmod\,2)$,
\begin{align*}
\sum_{\sigma\in A_j}(-1)^{\ell(\sigma)}x^{L(\sigma)}+\sum_{\sigma\in A_{j+1}}(-1)^{\ell(\sigma)}x^{L(\sigma)}
=\sum_{\sigma\in A_j}\left((-1)^{\ell(\sigma)}x^{L(\sigma)}+(-1)^{\ell(\sigma'')}x^{L(\sigma'')}\right)=0,
\end{align*}
where $\sigma''=(\sigma(i+j),\sigma(i+2k+2))(\overline{\sigma(i+j)},\overline{\sigma(i+2k+2)})\sigma$, since $\ell(\sigma'')=\ell(\sigma)-1$ and $L(\sigma'')=L(\sigma)$.
Therefore,
\begin{align*}
\sum_{\substack{\sigma\in \mathcal{D}_n^I,\, \sigma(b)=n'}}(-1)^{\ell(\sigma)}x^{L(\sigma)}
=\sum_{\substack{\sigma\in \mathcal{D}_n^I,\, \sigma(b)=n'\\ \sigma(i+2k+1)<\sigma(i+2k+2)}}(-1)^{\ell(\sigma)}x^{L(\sigma)}
=\sum_{\substack{\sigma\in \mathcal{D}_n^{I\cup \widetilde I},\, \sigma(b)=n'}}(-1)^{\ell(\sigma)}x^{L(\sigma)}.
\end{align*}
The second equality can be proved analogously.
\end{proof}

The following result reveals the effect of shifting the connected component of $I$ that contains $0$, on $\mathcal{D}_n^I(x)$.
\begin{lemma}\label{lem:a0=even}
Let $I\subseteq[0,n-1]$ with $I_0=[0,a_0-1]$, where $a_0\in[2, n-1]$.
If $a_0\equiv 0 (\bmod\,2)$ and $a_0+1\not\in I$, then
\begin{enumerate}
\item \label{itemeq:DnxIa0a}
$\mathcal{D}_n^I(x)=(1+x^{a_0})\mathcal{D}_n^{I\cup \{a_0\}}(x)$,\vspace{2mm}
\item \label{itemeq:DnxIa0=nb} for any  $b\in[a_0+2,n]$ and $n'\in \{n,\overline{n}\}$, we have
$$
\sum_{\substack{\sigma\in \mathcal{D}_n^I,\, \sigma(b)=n'}}(-1)^{\ell(\sigma)}x^{L(\sigma)}
=(1+x^{a_0})\sum_{\substack{\sigma\in \mathcal{D}_n^{I\cup \{a_0\}},\, \sigma(b)=n'}}(-1)^{\ell(\sigma)}x^{L(\sigma)}.
$$
\end{enumerate}
\end{lemma}
\begin{proof}
\eqref{itemeq:DnxIa0a}
Since $I_0=[0,a_0-1]$ and  $a_0\in[2,n-1]$, we see that $|\sigma(1)|<\sigma(2)<\cdots<\sigma(a_0)$ for any $\sigma\in \mathcal{D}_n^I$.
Let $A_j=\{\sigma\in \mathcal{D}_n^I\mid \sigma(j-1)<|\sigma(a_0+1)|<\sigma(j)\}$ for $j\in[3,a_0]$. Then
\begin{align}\label{eq:pfdnix4terms}
\mathcal{D}_n^I(x)= &\sum_{\substack{\sigma\in \mathcal{D}_n^I\\ |\sigma(a_0+1)|<|\sigma(1)|}}(-1)^{\ell(\sigma)}x^{L(\sigma)}+\sum_{\substack{\sigma\in \mathcal{D}_n^I\\ |\sigma(1)|<|\sigma(a_0+1)|<\sigma(2)}}(-1)^{\ell(\sigma)}x^{L(\sigma)}\notag\\
&\hspace{3mm}+\sum_{j=3}^{a_0}\sum_{\sigma\in A_j}(-1)^{\ell(\sigma)}x^{L(\sigma)}
+\sum_{\substack{\sigma\in \mathcal{D}_n^I\\ \sigma(a_0)<|\sigma(a_0+1)|}}(-1)^{\ell(\sigma)}x^{L(\sigma)}.
\end{align}
For each $\sigma\in \mathcal{D}_n^I$ with $|\sigma(a_0+1)|<|\sigma(1)|$, define
$$\sigma'=(\barr{|\sigma(1)|},\barr{|\sigma(a_0+1)|})(|\sigma(1)|,|\sigma(a_0+1)|)\sigma.$$
Since $a_0+1\not\in I$, the map $\sigma\mapsto \sigma'$ is a bijection from $\{\sigma\in \mathcal{D}_n^I\mid |\sigma(a_0+1)|<|\sigma(1)|\}$
to $\{\sigma\in \mathcal{D}_n^I\mid |\sigma(1)|<|\sigma(a_0+1)|<\sigma(2)\}$.
But $a_0\equiv 0 (\bmod\,2)$, a straightforward computation shows that $\ell(\sigma)=\ell(\sigma')\pm 1$ and $L(\sigma)=L(\sigma')$, so
the sum of the first and second terms on the right-hand side of \eqref{eq:pfdnix4terms} is
\begin{align*}
\sum_{\substack{\sigma\in \mathcal{D}_n^I\\ |\sigma(a_0+1)|<|\sigma(1)|}}
   \left((-1)^{\ell(\sigma)}x^{L(\sigma)}+(-1)^{\ell(\sigma')}x^{L(\sigma')}\right)
=0.
\end{align*}
Similarly, for any $j\in[3,a_0]$ with $j\equiv 1(\bmod\,2)$, we have
\begin{align*}
\sum_{\sigma\in A_j}(-1)^{\ell(\sigma)}x^{L(\sigma)}
+\sum_{\sigma\in A_{j+1}}(-1)^{\ell(\sigma)}x^{L(\sigma)}
=\sum_{\sigma\in A_j}\left((-1)^{\ell(\sigma)}x^{L(\sigma)}
+(-1)^{\ell(\sigma'')}x^{L(\sigma'')}\right)
=0,
\end{align*}
where  $\sigma''=(\barr{|\sigma(j)|},\barr{|\sigma(a_0+1)|})(|\sigma(j)|,|\sigma(a_0+1)|)\sigma$.
Since $a_0\equiv 0 (\bmod\,2)$, we see that the third summand on the right-hand side of \eqref{eq:pfdnix4terms} is $0$.
Notice that $\sigma(a_0)>0$.
Putting these together with \eqref{eq:pfdnix4terms} gives
\begin{align*}
\mathcal{D}_n^I(x)=& \sum_{\substack{\sigma\in \mathcal{D}_n^I\\ \sigma(a_0)<|\sigma(a_0+1)|}}(-1)^{\ell(\sigma)}x^{L(\sigma)}\\
=& \sum_{\substack{\sigma\in \mathcal{D}_n^I\\ \sigma(a_0)<\sigma(a_0+1)}}(-1)^{\ell(\sigma)}x^{L(\sigma)}+\sum_{\substack{\sigma\in \mathcal{D}_n^I\\ \sigma(a_0+1)<\barr{\sigma(a_0)}}}(-1)^{\ell(\sigma)}x^{L(\sigma)}\\
=& \sum_{\substack{\sigma\in \mathcal{D}_n^I\\ \sigma(a_0)<\sigma(a_0+1)}}\left((-1)^{\ell(\sigma)}x^{L(\sigma)}+(-1)^{\ell(\sigma''')}x^{L(\sigma''')}\right),
\end{align*}
where $\sigma'''=(1,\barr1)(\sigma(a_0+1),\barr{\sigma(a_0+1)})\sigma$. Since $\sigma(a_0+1)>\sigma(a_0)>0$,
it follows from Lemma \ref{lem:ellLsigma'a=n} that $\ell(\sigma''')=\ell(\sigma)+2a_0$ and $L(\sigma''')=L(\sigma)+a_0$. Therefore,
\begin{align*}
\mathcal{D}_n^I(x)=(1+x^{a_0})\sum_{\substack{\sigma\in \mathcal{D}_n^I\\ \sigma(a_0)<\sigma(a_0+1)}}(-1)^{\ell(\sigma)}x^{L(\sigma)}
=(1+x^{a_0}) \mathcal{D}_n^{I\cup\{a_0\}}(x),
\end{align*}
as desired.

\eqref{itemeq:DnxIa0=nb} Completely analogous to the proof of part \eqref{itemeq:DnxIa0a}, where the condition $\sigma\in \mathcal{D}_n^I$ is replaced by $\sigma\in \mathcal{D}_n^I$ and
$\sigma(b)=n'$.
\end{proof}

\subsection{Compressing}
Lemma \ref{prop:shift[0,n-1]}\eqref{eq:prop3.7dnx0} indicates that some connected components of $I$ can be compressed and shifted to the left or to the right in $[0,n-1]$,
without changing the polynomial $\mathcal{D}_n^I(x)$.
In this subsection, we push further this idea by introducing the following useful notion of compressed subsets.
\begin{defn}\label{defn:compressedset}
Let $a_0<a_1<\cdots<a_s$ be integers in $[0,n]$. The subset
$$[0,a_0-1]\cup[a_0+1,a_1-1]\cup \cdots\cup [a_{s-1}+1,a_s-1]$$
 of $[0,n-1]$ is said to be
\emph{compressed} if $a_0\equiv a_1\equiv \cdots\equiv a_s(\bmod\ 2)$, that is, $a_0,a_1,\ldots,a_s$ have the same parity.
 \end{defn}

Given a subset $I$ of $[0,n-1]$ with connected components $I_0,I_1,\ldots,I_s$, we associate it with the following compressed subset, denoted $C(I)$, defined by
$$
C(I):=[0,b_0-1]\cup[b_0+1,b_1-1]\cup \cdots\cup [b_{s-1}+1,b_s-1],
$$
where $b_0=|I_0|$ and
$$
b_{k}=|I_0|+\sum_{i=1}^k2\left\lfloor\frac{|I_i|+1}{2}\right\rfloor \quad \text{for all} \quad k\in [s].
$$
For example, let $n=20$. If
$$I'=[0,4]\cup[6,9]\cup[11,15]\cup[17,19]\quad \text{and}\quad I''=[2,3]\cup[6,9]\cup[11,15]\cup[17,19],$$
then
$$
C(I')=[0,4]\cup[6,8]\cup[10,14]\cup[16,18]\quad \text{and}\quad C(I'')=[1]\cup[3,5]\cup[7,11]\cup[13,15].
$$
Let $J=C(I)$, and let $J_0=[0,b_0-1]$ and $J_k=[b_{k-1}+1,b_k-1]$ for all $k\in[s]$. It is straightforward to check that
\begin{align*}
|J_k|=b_k-b_{k-1}-1=2\left\lfloor\frac{|I_k|+1}{2}\right\rfloor-1
 =\begin{cases}
	|I_k|, &{\rm if}\ |I_k|\equiv1(\bmod\,2),\\
	|I_k|-1, &{\rm if}\ |I_k|\equiv0(\bmod\,2).
\end{cases}
\end{align*}
Thus,
$\left\lfloor\frac{|J_k|+1}{2}\right\rfloor=\left\lfloor\frac{|I_k|+1}{2}\right\rfloor$ for all $k\in[0,s]$, and hence
\begin{align}\label{eq:compressm}
    m_I=m_J
    =\begin{cases}
    \frac{b_s+1}{2}, &{\rm if}\ |I_0|\equiv1(\bmod\,2),\\
    \frac{b_s}{2}, &{\rm if}\ |I_0|\equiv0(\bmod\,2).
    \end{cases}
\end{align}

Notice that the previous definition of a compressed subset of $[0,n-1]$ is similar to, but slightly different from, the one given in \cite[Section 4]{BC17T},
where a compressed set is only defined for a subset of $[n-1]$.

The following lemma suggests that if $|I_0|\geq2$, then the computation of $\mathcal{D}_n^I(x)$ can be reduced to the case where $I$ is a compressed subset of $[0,n-1]$.

\begin{lemma}\label{lem:tocompressset}
Let $I\subseteq[0,n-1]$ with $I_0=[0,a_0-1]$, where $a_0\geq2$. Then
\begin{enumerate}
\item\label{item:comprs1}  $\mathcal{D}_n^I(x)=\mathcal{D}_n^{C(I)}(x)$.
\vspace{3mm}
\item\label{item:comprs2}
$
\sum\limits_{\substack{\sigma\in \mathcal{D}_n^I,\, \sigma(a_0)=n}}(-1)^{\ell(\sigma)}x^{L(\sigma)}=
\sum\limits_{\substack{\sigma\in \mathcal{D}_n^{C(I)},\,\sigma(a_0)=n}}(-1)^{\ell(\sigma)}x^{L(\sigma)}.
$
\end{enumerate}
\end{lemma}
\begin{proof}
\eqref{item:comprs1}
It is routine to check that $C(I)=I$ if $I$ is a compressed subset. So it suffices to assume that $I$ is not compressed. Let $I_0$ and $I_k=[a_{k}+1,b_k-1]$, $k\in[s]$, are all connected components of $I$, where $b_{k-1}< a_{k}+1\leq b_{k}-1$, with the notation $b_0=a_0$.

If there exists $k\in[s]$ such that $|I_k|\equiv0(\bmod\ 2)$, then put $I'=I\setminus\{b_{k}-1\}$. By Lemma \ref{prop:shift[0,n-1]}\eqref{eq:prop3.7dnx0}, we see that $\mathcal{D}_n^I(x)=\mathcal{D}_n^{I'}(x)$. But now $|I_k'|=|I_k|-1$ is odd, so we have $\left\lfloor\frac{|I_k'|+1}{2}\right\rfloor=\left\lfloor\frac{|I_k|+1}{2}\right\rfloor$.
Thus, we may assume that $|I_k|\equiv1(\bmod\ 2)$ for all $k\in[s]$.

Since $I$ is not compressed, there exists an integer $j\in[s]$, such that $b_{j-1}<a_{j}$.
Suppose that $j$ is chosen to be minimal with this property.
Then $a_{j}-1$, $a_{j}$ and $b_{j}$ are not elements of $I$.
Let $I''=(I\setminus\{b_j-1\})\cup\{a_{j}\}$.
Since $|I_j|\equiv1(\bmod\ 2)$, we have $a_j\equiv b_j(\bmod\, 2)$,
it follows from Lemma \ref{prop:shift[0,n-1]}\eqref{eq:prop3.7dnx0} that
$\mathcal{D}_n^I(x)=\mathcal{D}_n^{I''}(x).$
Continuing in this way, we can get a compressed subset $J\subseteq[0,n-1]$ such that $\mathcal{D}_n^I(x)=\mathcal{D}_n^{J}(x)$.
Moreover, by the construction of $J$, we see that $J_0=I_0$ and
\begin{align*}
|J_k|=\begin{cases}
|I_k|-1, &{\rm if}\ |I_k|\equiv0(\bmod\,2),\\
|I_k|, &{\rm if}\ |I_k|\equiv1(\bmod\,2)
\end{cases}
\end{align*}
for all $k\in[s]$,
so that $|J_k|=2\left\lfloor\frac{|I_k|+1}{2}\right\rfloor-1$. Hence $C(I)=J$, and the proof follows.

\eqref{item:comprs2}
Since if $|I_0|\geq2$, then $I_0=[0,a_0-1]$ is a common connected component of $I$ and $C(I)$, the proof of part \eqref{item:comprs1} is still valid if we substitute the condition
$\sigma\in \mathcal{D}_n^{I}$ by
the condition $\sigma\in \mathcal{D}_n^{I}$ with $\sigma(a_0)=n$ and apply Lemma \ref{prop:shift[0,n-1]}\eqref{eq:prop3.7dnx0a=n}.
\end{proof}

The proof of the next result, which will be used later, is analogous to that of Lemma \ref{lem:a0=even}.

\begin{lemma}\label{lem:a0=oddn}
Let $I\subseteq[0,n-2]$ be compressed with $I_0=[0,a_0-1]$, where $a_0\in[3,n-1]$. If $a_0\equiv 1 (\bmod\,2)$, then
\begin{align*}
\sum_{\substack{\sigma\in \mathcal{D}_n^I,\, \sigma(a_0)=n}}(-1)^{\ell(\sigma)}x^{L(\sigma)}=
(-1)^{n-1}2x^{\left\lfloor\frac{n}{2}\right\rfloor}\mathcal{D}_{n-1}^I(x).
\end{align*}
\end{lemma}
\begin{proof}
Let $I=[0,a_0-1]\cup[a_0+1,a_1-1]\cup\cdots\cup[a_{s-1}+1,a_s-1]$, where $a_s\leq n-1$. It follows from Lemma
\ref{prop:shift[0,n-1]}\eqref{eq:prop3.7dnx0a=n} that
\begin{align*}
\sum_{\substack{\sigma\in \mathcal{D}_n^I,\, \sigma(a_0)=n}}(-1)^{\ell(\sigma)}x^{L(\sigma)}
%=\sum_{\substack{\sigma\in \mathcal{D}_n^{C(I)},\, \sigma(a_0)=n}}(-1)^{\ell(\sigma)}x^{L(\sigma)}
=\sum_{\substack{\sigma\in \mathcal{D}_n^{\overline I},\, \sigma(a_0)=n}}(-1)^{\ell(\sigma)}x^{L(\sigma)},
	\end{align*}
where $\overline I=[0,a_0-1]\cup[a_0+1,a_1-1]\cup\cdots\cup[a_{s-2}+1,a_{s-1}-1]\cup[a_{s-1}+2,a_s]$.
Repeated applications of Lemma \ref{prop:shift[0,n-1]}\eqref{eq:prop3.7dnx0a=n} give that
	\begin{align}\label{nsigmaa0=nto-1}
		\sum_{\substack{\sigma\in \mathcal{D}_n^I,\, \sigma(a_0)=n}}(-1)^{\ell(\sigma)}x^{L(\sigma)}
=\sum_{\substack{\sigma\in \mathcal{D}_n^{\widetilde I},\, \sigma(a_0)=n}}(-1)^{\ell(\sigma)}x^{L(\sigma)},
	\end{align}
where $\widetilde I=[0,a_0-1]\cup[a_0+2,a_1]\cup\cdots\cup[a_{s-1}+2,a_s]$.
For any $\sigma\in \mathcal{D}_n^{\widetilde I}$, we have
$|\sigma(1)|<\sigma(2)<\cdots<\sigma(a_0)$, and hence $|\sigma(1)|<|\sigma(2)|<\cdots<|\sigma(a_0)|$.
Then
\begin{align*}
\sum_{{\sigma\in \mathcal{D}_n^{\widetilde I},\, \sigma(a_0)=n}}(-1)^{\ell(\sigma)}x^{L(\sigma)}
=&\,\sum_{\substack{\sigma\in \mathcal{D}_n^{\widetilde I},\, \sigma(a_0)=n\\ |\sigma(a_0+1)|<|\sigma(1)|}}(-1)^{\ell(\sigma)}x^{L(\sigma)}
+\sum_{j=2}^{a_0}\sum_{\sigma\in A_j}(-1)^{\ell(\sigma)}x^{L(\sigma)},
\end{align*}
where $A_j=\{\sigma\in \mathcal{D}_n^{\widetilde I},\mid \sigma(a_0)=n,\, |\sigma(j-1)|<|\sigma(a_0+1)|<|\sigma(j)|\}$.
Since $a_0\equiv 1 (\bmod\,2)$ and $a_0\in[3, n-1]$, we see that for any $j\in[2,a_0]$ with $j\equiv 0(\bmod\,2)$,
\begin{align*}%\label{eq:|sigma(2j-1)|<|sigma(a_0+1)|<|sigma(2j)|11}
&\,\sum_{\sigma\in A_j} (-1)^{\ell(\sigma)}x^{L(\sigma)}
+\sum_{\sigma\in A_{j+1}}(-1)^{\ell(\sigma)}x^{L(\sigma)}
=\sum_{\sigma\in A_j}
   \left((-1)^{\ell(\sigma)}x^{L(\sigma)}+(-1)^{\ell(\sigma')}x^{L(\sigma')}\right)=0,
\end{align*}
where $\sigma'=(\barr{|\sigma(j)|},\barr{|\sigma(a_0+1)|})(|\sigma(j)|,|\sigma(a_0+1)|)\sigma$.
Thus,
\begin{align*}
\sum_{\substack{\sigma\in \mathcal{D}_n^{\widetilde I},\, \sigma(a_0)=n}}(-1)^{\ell(\sigma)}x^{L(\sigma)}=&
  \sum_{\substack{\sigma\in \mathcal{D}_n^{\widetilde I},\, \sigma(a_0)=n\\ |\sigma(a_0+1)|<|\sigma(1)|}}(-1)^{\ell(\sigma)}x^{L(\sigma)}\\
=& \sum_{\substack{\sigma\in \mathcal{D}_n^{\widetilde I},\, \sigma(a_0)=n\\ |\sigma(a_0+1)|<\sigma(1)}}(-1)^{\ell(\sigma)}x^{L(\sigma)}
+ \sum_{\substack{\sigma\in \mathcal{D}_n^{\widetilde I},\, \sigma(a_0)=n\\ |\sigma(a_0+1)|<\barr{\sigma(1)}}}(-1)^{\ell(\sigma)}x^{L(\sigma)}\\
=& \sum_{\sigma\in \mathcal{D}_{n-1}^{I}}\left((-1)^{\ell(\sigma'')}x^{L(\sigma'')}+(-1)^{\ell(\sigma''')}x^{L(\sigma''')}\right),
\end{align*}
where
\begin{align*}
\sigma''=&\ \sigma(2)\sigma(3)\cdots\sigma(a_0)\hspace{1pt}n\hspace{1pt}\sigma(1)\sigma(a_0+1)\sigma(a_0+2)\cdots\sigma(n-1),\\
%\intertext{and}
\sigma'''=&\ \barr{\sigma(2)}\sigma(3)\cdots\sigma(a_0)\hspace{1pt}n\hspace{1pt}\barr{\sigma(1)}\sigma(a_0+1)\sigma(a_0+2)\cdots\sigma(n-1).
\end{align*}
For any $\sigma\in \mathcal{D}_{n-1}^{I}$, it follows from $|\sigma(1)|<\sigma(2)<\cdots<\sigma(a_0)$ that
\delete{
\begin{align*}
\Inv(\sigma)=\ &\{(1,j)\,|\,j\in[a_0+1, n-1],\sigma(1)>\sigma(j)\}\\
&\biguplus\{(i,j)\,|\,i\in[2,a_0],j\in[a_0+1, n-1], \sigma(i)>\sigma(j)\}\\
&\biguplus\{(a_0+i,a_0+j)\,|\,(i,j)\in\Inv(\sigma(a_0+1)\cdots\sigma(n-1))\},\\
\Nsp(\sigma)=\ &\{(1,j)\,|\,j\in[a_0+1, n-1],\sigma(1)+\sigma(j)<0\}\\
&\biguplus\{(i,j)\,|\,i\in[2, a_0],j\in[a_0+1, n-1], \sigma(i)+\sigma(j)<0\}\\
&\biguplus\{(a_0+i,a_0+j)\,|\,(i,j)\in\Nsp(\sigma(a_0+1)\cdots\sigma(n-1))\},\\
\Inv(\sigma')=\ &\{(i,a_0+1)\,|\,i\in[a_0]\}\uplus\{(a_0,j)\,|\,j\in[a_0+2, n]\}\\
&\biguplus\{(a_0+1,j)\,|\,j\in[a_0+2, n],\sigma(1)>\sigma(j-1)\}\\
&\biguplus\{(i,j)\,|\,i\in[a_0-1],j\in[a_0+2,n], \sigma(i+1)>\sigma(j-1)\}\\
&\biguplus\{(a_0+1+i,a_0+1+j)\,|\,(i,j)\in\Inv(\sigma(a_0+1)\cdots\sigma(n-1))\}
\intertext{and}
\Nsp(\sigma')=\ &\{(a_0+1,j)\,|\,j\in[a_0+2, n],\sigma(1)+\sigma(j-1)<0\}\\
&\biguplus\{(i,j)\,|\,i\in[a_0-1],j\in[a_0+2, n], \sigma(i+1)+\sigma(j-1)<0\}\\
&\biguplus\{(a_0+1+i,a_0+1+j)\,|\,(i,j)\in\Nsp(\sigma(a_0+1)\cdots\sigma(n-1))\}.
\end{align*}
Thus, }
$\ell(\sigma'')-\ell(\sigma)=n-1$. Since $a_0\equiv 1(\bmod\ 2)$, we have
\begin{align*}
L(\sigma'')-L(\sigma)=&\oinv(\sigma'')-\oinv(\sigma)\\
=&\ |\{i\in[a_0]\,|\,i\not\equiv a_0+1(\bmod\ 2)\}|+|\{j\in[a_0+2, n]\,|\,j\not\equiv a_0(\bmod\ 2)\}|\\
%=&\ \frac{a_0+1}{2}+\left\lfloor\frac{n-a_0-1}{2}\right\rfloor\\
=&\,\left\lfloor\frac{n}{2}\right\rfloor.
\end{align*}
In an analogous manner, we see that $\ell(\sigma''')-\ell(\sigma)=n-1$ and $L(\sigma''')-L(\sigma)=\left\lfloor\frac{n}{2}\right\rfloor$.
Thus,
	\begin{align*}
\sum_{\substack{\sigma\in \mathcal{D}_n^{\widetilde I},\, \sigma(a_0)=n}}(-1)^{\ell(\sigma)}x^{L(\sigma)}=(-1)^{n-1}2x^{\left\lfloor\frac{n}{2}\right\rfloor}\mathcal{D}_{n-1}^{I}(x),
	\end{align*}
and the proof follows from \eqref{nsigmaa0=nto-1}.
\end{proof}
\delete{
\begin{lemma}\label{lem:sum=000}
Let $n\in\mathbb{P}$, and $I=[0,a_0-1]\cup [a_{0}+1,a_1-1]\cup\cdots\cup [a_{s-1}+1,a_s-1]$ be a compressed subset of $[0,n-2]$, where $\{a_0,a_1,\ldots,a_s\}_{<}\subseteq [n]$ and $a_0\geq3$ is an odd positive integer. Then
\begin{align*}%\label{eq:isigmaa0=n}
\sum_{\substack{\sigma\in \mathcal{D}_n^I,\, \sigma(a_i+1)=\overline{n}}}(-1)^{\ell(\sigma)}x^{L(\sigma)}
=0\quad \text{for all}\ i\in[0,s-1].
\end{align*}
\end{lemma}
\begin{proof}
It follows from Lemma \ref{lem:a0=even}\eqref{itemeq:DnxIa0=nb} that for each $i\in[0,s-1]$, we have
\begin{align*}%\label{eq:isigmaa0=n}
\sum_{\substack{\sigma\in \mathcal{D}_n^I,\, \sigma(a_i+1)=\overline{n}}}(-1)^{\ell(\sigma)}x^{L(\sigma)}
&=\frac{1}{1+x^{a_0-1}}\sum_{\substack{\sigma\in \mathcal{D}_n^{I\backslash\{a_0-1\}},\, \sigma(a_i+1)=\overline{n}}}(-1)^{\ell(\sigma)}x^{L(\sigma)}
&&{\rm  by\ Lemma\ \ref{lem:a0=even}\eqref{itemeq:DnxIa0=nb}}\\
&=\frac{1}{1+x^{a_0-1}}\sum_{\substack{\sigma\in \mathcal{D}_n^{\widetilde{I}},\, \sigma(a_i+1)=\overline{n}}}(-1)^{\ell(\sigma)}x^{L(\sigma)}&&{\rm  by\ Lemma\ \ref{prop:shift[0,n-1]}\eqref{eq:prop3.7dnx0a=n}}\\
&=0&&{\rm  by\ Lemma\ \ref{lem:ai=oddn}\eqref{item:sigmaa=-na}},
\end{align*}
where $\widetilde{I}=[0,a_0-2]\cup [a_{0},a_1-2]\cup\cdots\cup[a_{i-1},a_i-2]\cup[a_{i}+1,a_{i+1}-1]\cup\cdots\cup [a_{s-1}+1,a_s-1]$.
\end{proof}
}

The following two results indicate that the computation of $\mathcal{D}_n^I(x)$ amounts to the computation of the sign-twisted generating functions over some elements for which $n$ (or $\overline{n}$) is in certain positions.
\begin{lemma}[\cite{BC17D}, Lemma 3.3]\label{lem:eliminatesigma(a)=n}
Let $n\geq 3$, $I\subseteq[0,n-1]$, and let $a\in[2,n-1]$ such that $a+1\notin I$. Suppose that the following hold: if $a=3$, then $0,1\notin I$; if $a\geq4$, then $a-2\notin I$. Then for each $n'\in\{n,\overline{n}\}$, we have
$
\sum\limits_{\substack{\sigma\in \mathcal{D}_n^I,\, \sigma(a)=n'}}(-1)^{\ell(\sigma)}x^{L(\sigma)}=0.
$
\end{lemma}

\begin{lemma}\label{lem:ai=oddn}
Let $I\subseteq[0,n-1]$, and let $[i,k]$ be a connected component of $I$, where $i\geq 1$.
\begin{enumerate}
\item\label{item:sigmaa=-na}
If $i\geq2$, then
\begin{align*}
\sum_{\substack{\sigma\in \mathcal{D}_n^I,\, \sigma(i)=\overline{n}}}(-1)^{\ell(\sigma)}x^{L(\sigma)}=
\begin{cases}
\sum\limits_{\substack{\sigma\in \mathcal{D}_n^I,\, \sigma(i)=\overline{n}\\ \sigma(k+1)<\sigma(i-1)}}(-1)^{\ell(\sigma)}x^{L(\sigma)}, &{\rm if}\ i\not\equiv k(\bmod\,2),\\
\hspace{8mm}0, & {\rm if}\ i\equiv k(\bmod\,2).
\end{cases}
\end{align*}
\item\label{lemeq:sigma(b)=n2} If $k\leq n-2$ and $k+2\not\in I$, then
\begin{align*}
\sum_{\substack{\sigma\in \mathcal{D}_n^I,\, \sigma(k+1)=n}}(-1)^{\ell(\sigma)}x^{L(\sigma)}=
\begin{cases}
\sum\limits_{\substack{\sigma\in \mathcal{D}_n^I,\, \sigma(k+1)=n\\ \sigma(k+2)<\sigma(i)}}(-1)^{\ell(\sigma)}x^{L(\sigma)}, &{\rm if}\ i\not\equiv k(\bmod\,2),\\
\hspace{8mm}0, & {\rm if}\ i\equiv k(\bmod\,2).
\end{cases}
\end{align*}
\end{enumerate}
\end{lemma}
\begin{proof}
\eqref{item:sigmaa=-na} Since $i\geq2$ and $[i,k]$ is a connected component of $I$, we have
$\sigma(i)<\sigma(i+1)<\cdots<\sigma(k+1)$ for any $\sigma\in \mathcal{D}_n^I$.
Let $A_j=\{\sigma\in \mathcal{D}_n^{I}\mid \sigma(i)=\overline{n},\sigma(i+j-1)<\sigma(i-1)<\sigma(i+j)\}$, where $j\in[k-i+1]$.
Then, by considering the value of $\sigma(i-1)$, we obtain that
\begin{align*}%\label{eq:sigma(a)=-n1}
\sum_{\sigma\in \mathcal{D}_n^{ I},\sigma(i)=\overline{n}}(-1)^{\ell(\sigma)}x^{L(\sigma)}
=\sum_{j=1}^{k-i+1}\sum_{\sigma\in A_j}(-1)^{\ell(\sigma)}x^{L(\sigma)}
+\sum_{\substack{\sigma\in \mathcal{D}_n^{I},\,\sigma(i)=\overline{n}\\ \sigma(k+1)<\sigma(i-1)}}(-1)^{\ell(\sigma)}x^{L(\sigma)}.
\end{align*}
For each $j\in[k-i+1]$ with $j\equiv 1(\bmod\,2)$, we  have
\begin{align*}
\sum_{\sigma\in A_j} (-1)^{\ell(\sigma)}x^{L(\sigma)}
+\sum_{\sigma\in A_{j+1}}(-1)^{\ell(\sigma)}x^{L(\sigma)}
=\sum_{\sigma\in A_j}\left( (-1)^{\ell(\sigma)}x^{L(\sigma)}
+(-1)^{\ell(\sigma')}x^{L(\sigma')}\right)
=0,
\end{align*}
since  $\ell(\sigma')=\ell(\sigma)+ 1$ and $L(\sigma')=L(\sigma)$, where  $\sigma'=\sigma(i-1,i+j)(\overline{i-1},\overline{i+j})$.
Thus, if $i\not\equiv k(\bmod\,2)$, then
\begin{align*}
\sum_{\substack{\sigma\in \mathcal{D}_n^I,\, \sigma(i)=\overline{n}}}(-1)^{\ell(\sigma)}x^{L(\sigma)}=
\sum\limits_{\substack{\sigma\in \mathcal{D}_n^I,\, \sigma(i)=\overline{n}\\ \sigma(k+1)<\sigma(i-1)}}(-1)^{\ell(\sigma)}x^{L(\sigma)}.
\end{align*}
If $i\equiv k(\bmod\,2)$, then
\begin{align*}
\sum_{\substack{\sigma\in \mathcal{D}_n^{ I},\, \sigma(i)=\overline{n}}}(-1)^{\ell(\sigma)}x^{L(\sigma)}
&=\sum_{\sigma\in A_{k-i+1}}(-1)^{\ell(\sigma)}x^{L(\sigma)}
+\sum_{\substack{\sigma\in \mathcal{D}_n^{I},\, \sigma(i)=\overline{n}\\ \sigma(k+1)<\sigma(i-1)}}(-1)^{\ell(\sigma)}x^{L(\sigma)}\\
&=\sum_{\sigma\in A_{k-i+1}}\left((-1)^{\ell(\sigma)}x^{L(\sigma)}+(-1)^{\ell(\sigma'')}x^{L(\sigma'')}\right)\\
&=0,
\end{align*}
since $\ell(\sigma'')=\ell(\sigma)+1$ and $L(\sigma'')=L(\sigma)$,
where $\sigma'=\sigma(i-1,k+1)(\barr{i-1},\barr{k+1})$.

\eqref{lemeq:sigma(b)=n2} If $k\leq n-2$ and $k+2\not\in I$, then
\begin{align*}%\label{eq:sigma(a)=-n2}
\sum_{\substack{\sigma\in \mathcal{D}_n^I\\ \sigma(k+1)=n}}(-1)^{\ell(\sigma)}x^{L(\sigma)}=\sum_{j=1}^{k-i+1}\sum_{\substack{\sigma\in \mathcal{D}_n^{I},\, \sigma(k+1)=n\\ \sigma(i+j-1)<\sigma(k+2)<\sigma(i+j)}}(-1)^{\ell(\sigma)}x^{L(\sigma)}
+\sum_{\substack{\sigma\in \mathcal{D}_n^{I},\, \sigma(k+1)=n\\ \sigma(k+2)<\sigma(i)}}(-1)^{\ell(\sigma)}x^{L(\sigma)}.
\end{align*}
The remainder of the proof is analogous to that of part \eqref{item:sigmaa=-na}; the details are omitted.
\end{proof}

\delete{
\begin{remark}
In the proof of Lemma \ref{lem:ai=oddn}, we use the approach of sign-reversing involutions, which is
one of the standard combinatorial techniques for eliminating cancellations.
Let $A$ be a set equipped with a sign function ${\rm sgn}:A\rightarrow\{1,-1\}$. An involution $\varphi:A\rightarrow A$ is \emph{sign reversing} if
${\rm sgn}\,\varphi(a)=-{\rm sgn}\, a $ when $\varphi(a)\neq a$. It follows that
$$
\sum_{a\in A}{\rm sgn}a=\sum_{a\in F}{\rm sgn}a,
$$
where $F$ is the set of fixed points of $\varphi$. In particular, if $\varphi$ has no fixed points, then $\sum\limits_{a\in A}{\rm sgn}a=0$.
Such a method of proof is called a \emph{sign-reversing involution argument}.
\end{remark}}

The main results of this section are the following two propositions, which will be used in a crucial way to deduce the formulas for $\mathcal{D}_n^I(x)$.

\begin{prop}\label{prop:a0=odd}
Let $I\subseteq[0,n-1]$ with $|I_0|\geq2$, and let $J=C(I)$. If $|I_0|\equiv 1(\bmod\ 2)$, then
\begin{align*}
\mathcal{D}_n^I(x)=\mathcal{D}_{2m_J-1}^{J}(x)\prod_{j=2m_J}^n\left(1+(-1)^{j-1}x^{\left\lfloor\frac{j}{2}\right\rfloor}\right)^2.
\end{align*}
\end{prop}
\begin{proof}
By Lemma \ref{lem:tocompressset}\eqref{item:comprs1}, $\mathcal{D}_n^I(x)=\mathcal{D}_n^J(x)$.
Since $C(J)=J$ and $m_I=m_J$, it suffices to assume that $I$
is compressed and prove that
\begin{align*}
\mathcal{D}_n^I(x)=\mathcal{D}_{2m_I-1}^{I}(x)\prod_{j=2m_I}^n\left(1+(-1)^{j-1}x^{\left\lfloor\frac{j}{2}\right\rfloor}\right)^2.
\end{align*}
%Recall that $m_I=\sum_{k=0}^s\left\lfloor\frac{|I_k|+1}{2}\right\rfloor$.
Since $|I_0|\equiv 1(\bmod\ 2)$  with $|I_0|\geq2$, we assume that
$I_0=[0,a_0-1]$ and $I_k=[a_{k-1}+1,a_k-1]$ for all $k\in[s]$, where $3\leq a_0<a_1<\cdots<a_s\leq n$ with
$a_0\equiv a_1\equiv \cdots\equiv a_s\equiv 1(\bmod\,2)$.
It follows from \eqref{eq:compressm} that
$a_s=2m_I-1$, and hence $I\subseteq[0,2m_I-2]$.
We proceed by induction on $n$,
where the base case $n=a_s$ is trivial. Now assume that the desired identity is true for $n-1\geq a_s$.
Since $a_0\geq 3$, there follows $|\sigma(1)|<\sigma(2)$, so that $\sigma(1)\not\in\{n,\overline{n}\}$ for any $\sigma\in \mathcal{D}_n^I$.
Then, considering the positions of $n$ and $\overline{n}$ in the one-line notation of $\sigma$, it follows from Lemma \ref{lem:eliminatesigma(a)=n} that
\begin{align}\label{eq:finalboss}
\mathcal{D}_n^I(x)=&\ \sum_{i=0}^s\sum_{\substack{\sigma\in \mathcal{D}_n^I,\, \sigma(a_i)=n}}(-1)^{\ell(\sigma)}x^{L(\sigma)}
+\sum_{i=0}^{s-1}\sum_{\substack{\sigma\in \mathcal{D}_n^I,\, \sigma(a_i+1)=\overline n}}(-1)^{\ell(\sigma)}x^{L(\sigma)}\\
&  +\sum_{\substack{\sigma\in \mathcal{D}_n^I,\, |\sigma(a_s+1)|=n}}(-1)^{\ell(\sigma)}x^{L(\sigma)}
 +\sum_{\substack{\sigma\in \mathcal{D}_n^I,\, |\sigma(n)|=n}}(-1)^{\ell(\sigma)}x^{L(\sigma)},\notag
\end{align}
where the third summand exists if and only if $a_s\leq n-2$.

Since $a_0\geq3$ and $a_0\equiv 1(\bmod\ 2)$, it follows from Lemma \ref{lem:a0=oddn} that
$$
\sum_{\substack{\sigma\in \mathcal{D}_n^I,\, \sigma(a_0)=n}}(-1)^{\ell(\sigma)}x^{L(\sigma)}
= (-1)^{n-1}2x^{\left\lfloor\frac{n}{2}\right\rfloor}\mathcal{D}_{n-1}^{I}(x).
$$
For each $i\in[s]$, we have $a_{i}\leq n-1$, it follows from Lemmas \ref{prop:shift[0,n-1]}\eqref{eq:prop3.7dnx0a=n} and \ref{lem:ai=oddn}\eqref{lemeq:sigma(b)=n2} that
\begin{align*}%\label{eq:isigmaa0=n}
\sum_{\substack{\sigma\in \mathcal{D}_n^I,\, \sigma(a_i)=n}}(-1)^{\ell(\sigma)}x^{L(\sigma)}
=\sum_{\substack{\sigma\in \mathcal{D}_n^{\overline{I}},\, \sigma(a_i)=n}}(-1)^{\ell(\sigma)}x^{L(\sigma)}
=0,
\end{align*}
where $\overline{I}=[0,a_0-1]\cup\cdots\cup[a_{i-1}+1,a_i-1]\cup[a_{i}+2,a_{i+1}]\cup\cdots\cup [a_{s-1}+2,a_s].$

For each $i\in[0,s-1]$, we have
\begin{align*}%\label{eq:isigmaa0=n3}
\sum_{\sigma\in \mathcal{D}_n^I, \sigma(a_i+1)=\overline{n}}(-1)^{\ell(\sigma)}x^{L(\sigma)}
&=\frac{1}{1+x^{a_0-1}}\sum_{{\sigma\in \mathcal{D}_n^{I\backslash\{a_0-1\}}, \sigma(a_i+1)=\overline{n}}}(-1)^{\ell(\sigma)}x^{L(\sigma)}
&&{\rm  by\ Lemma\ \ref{lem:a0=even}\eqref{itemeq:DnxIa0=nb}}\\
&=\frac{1}{1+x^{a_0-1}}\sum_{{\sigma\in \mathcal{D}_n^{\widetilde{I}}, \sigma(a_i+1)=\overline{n}}}(-1)^{\ell(\sigma)}x^{L(\sigma)}&&{\rm  by\ Lemma\ \ref{prop:shift[0,n-1]}\eqref{eq:prop3.7dnx0a=n}}\notag\\
&=0,&&{\rm  by\ Lemma\ \ref{lem:ai=oddn}\eqref{item:sigmaa=-na}}\notag
\end{align*}
where $\widetilde{I}=[0,a_0-2]\cup [a_{0},a_1-2]\cup\cdots\cup[a_{i-1},a_i-2]\cup[a_{i}+1,a_{i+1}-1]\cup\cdots\cup [a_{s-1}+1,a_s-1]$.

If $a_s\leq n-2$, then the third term on the right-hand side of \eqref{eq:finalboss} exists, and it follows from Lemma\ \ref{lem:ellLsigma'a=n} that
\begin{align*}
\sum_{\substack{\sigma\in \mathcal{D}_n^I,\, |\sigma(a_s+1)|=n}}(-1)^{\ell(\sigma)}x^{L(\sigma)}
=&\,\sum_{\substack{\sigma\in \mathcal{D}_n^I,\, \sigma(a_s+1)=n}}(-1)^{\ell(\sigma)}x^{L(\sigma)}
+\sum_{\substack{\sigma\in \mathcal{D}_n^I,\, \sigma(a_s+1)=\overline{n}}}(-1)^{\ell(\sigma)}x^{L(\sigma)}\\
=&\,\sum_{\substack{\sigma\in \mathcal{D}_n^I,\, \sigma(a_s+1)=n}}\left((-1)^{\ell(\sigma)}x^{L(\sigma)}+(-1)^{\ell(\sigma')}x^{L(\sigma')}\right)\\
=&\,(1+x^{a_s+1})\sum_{\substack{\sigma\in \mathcal{D}_n^{I},\, \sigma(a_s+1)=n}}(-1)^{\ell(\sigma)}x^{L(\sigma)},
\end{align*}
since $\ell(\sigma')-\ell(\sigma)=2a_s$ and $L(\sigma')-L(\sigma)=a_s+1$, where $\sigma'=(1,\overline{1})(n,\overline{n})\sigma$. But
\begin{align*}
\sum_{\substack{\sigma\in \mathcal{D}_n^{I},\, \sigma(a_s+1)=n}}(-1)^{\ell(\sigma)}x^{L(\sigma)}
=&\,\frac{1}{1+x^{a_0-1}}\sum_{\substack{\sigma\in \mathcal{D}_n^{I\backslash\{a_{0}-1\}}\\ \sigma(a_s+1)=n}}(-1)^{\ell(\sigma)}x^{L(\sigma)}&&{\rm by\ Lemma\ \ref{lem:a0=even}\eqref{itemeq:DnxIa0=nb}}\\
=&\,\frac{1}{1+x^{a_0-1}}\sum_{\sigma\in \mathcal{D}_n^{J}, \sigma(a_s+1)=n}(-1)^{\ell(\sigma)}x^{L(\sigma)}
&&{\rm by\ Lemma\ \ref{prop:shift[0,n-1]}\eqref{eq:prop3.7dnx0a=n}}\\
=&\,0,&&{\rm by\ Lemma\ \ref{lem:eliminatesigma(a)=n}}
\end{align*}
where $J=[0,a_0-2]\cup[a_0,a_1-2]\cup\cdots\cup[a_{s-1},a_{s}-2].$
Hence the third term on the right-hand side of \eqref{eq:finalboss} is zero.

Notice that
\begin{align*}
\sum_{\substack{\sigma\in \mathcal{D}_n^I,\, |\sigma(n)|=n}}(-1)^{\ell(\sigma)}x^{L(\sigma)}
=&\,\sum_{\substack{\sigma\in \mathcal{D}_n^I,\, \sigma(n)=n}}(-1)^{\ell(\sigma)}x^{L(\sigma)}
  +\sum_{\substack{\sigma\in \mathcal{D}_n^I,\, \sigma(n)=\bar n}}(-1)^{\ell(\sigma)}x^{L(\sigma)}\\
=&\,\sum_{\substack{\sigma\in \mathcal{D}_n^I,\, \sigma(n)=n}}\left((-1)^{\ell(\sigma)}x^{L(\sigma)}+(-1)^{\ell(\sigma')}x^{L(\sigma')}\right),
\end{align*}
where $\sigma'=(1,\overline{1})(n,\overline{n})\sigma$. Then, by Lemma \ref{lem:ellLsigma'a=n},
\begin{align*}
\sum_{\substack{\sigma\in \mathcal{D}_n^I,\, |\sigma(n)|=n}}(-1)^{\ell(\sigma)}x^{L(\sigma)}
=\left(1+x^{2\left\lfloor\frac{n}{2}\right\rfloor}\right)\sum_{\substack{\sigma\in \mathcal{D}_n^I,\, \sigma(n)=n}}(-1)^{\ell(\sigma)}x^{L(\sigma)}
=\left(1+x^{2\left\lfloor\frac{n}{2}\right\rfloor}\right)\mathcal{D}_{n-1}^I(x).
\end{align*}
Combining with equation \eqref{eq:finalboss} gives
\begin{align*}
\mathcal{D}_n^I(x)
=& \left(1+(-1)^{n-1}x^{\left\lfloor\frac{n}{2}\right\rfloor}\right)^2\mathcal{D}_{n-1}^I(x)
= \mathcal{D}_{a_s}^I(x)\prod_{j=a_s+1}^n\left(1+(-1)^{j-1}x^{\left\lfloor\frac{j}{2}\right\rfloor}\right)^2,
\end{align*}
and the proof follows by induction and noting that $a_s=2m_I-1$.
\end{proof}

\begin{prop}\label{prop:a0=even}
Let  $I\subseteq[0,n-1]$ with $|I_0|\geq2$, and let $C(I)=[0,a_0-1]\cup[a_0+1,a_1-1]\cup\cdots\cup[a_{s-1}+1,a_s-1]$ with $a_s\leq n-1$. If $|I_0|\equiv 0(\bmod\ 2)$, then
\begin{align*}
\mathcal{D}_n^I(x)=(1+x^{a_0})\mathcal{D}_{2m_J-1}^{J}(x)\prod_{j=2m_J}^n\left(1+(-1)^{j-1}x^{\left\lfloor\frac{j}{2}\right\rfloor}\right)^2,
\end{align*}
where $J=[0,a_0]\cup[a_0+2,a_1]\cup\cdots\cup[a_{s-1}+2,a_s]$.
\end{prop}
\begin{proof}
Since $a_s\leq n-1$, it follows from Lemmas \ref{lem:tocompressset}\eqref{item:comprs1}, \ref{prop:shift[0,n-1]}\eqref{eq:prop3.7dnx0}, and \ref{lem:a0=even}\eqref{itemeq:DnxIa0a} that
$$\mathcal{D}_n^I(x)=\mathcal{D}_{n}^{C(I)}(x)=\mathcal{D}_{n}^{J\backslash\{a_0\}}(x)=(1+x^{a_0})\mathcal{D}_{n}^{J}(x).$$
But $J=C(J)$ and $|J_0|\equiv 1(\bmod\ 2)$, so
the proof now follows from Proposition \ref{prop:a0=odd}.
\end{proof}

\section{Chessboard elements and odd sandwiches}\label{sec:Ceos}

By \eqref{eq:compressm} and Propositions \ref{prop:a0=odd} and \ref{prop:a0=even}, in order to determine $\mathcal{D}_n^I(x)$ when $|I_0|\geq 2$, it suffices to compute $\mathcal{D}_n^J(x)$ for the case where $J$ is compressed with the form $J=[0,a_0-1]\cup[a_0+1,a_1-1]\cup\cdots\cup[a_{s-1}+1,n-1]$.
The task of this section is to show that $\mathcal{D}_n^J(x)$ can be obtained by computing $\mathcal{D}_n^{[0,a_0-1]\cup[a_0+1,n-1]}(x)$ (see Proposition \ref{lem:DnI=DnJxSn-aKeq}).

Following \cite{KV09,SV13}, we introduce the notions of chessboard elements of $\mathcal{D}_n$ and odd sandwiches of even signed permutations
that will play important roles in the paper.
We show in subsection \ref{subsec:Chessboardelements} that $\mathcal{D}_n^I(x)$ is supported on the set of chessboard elements of $\mathcal{D}_n^I$.
In subsection \ref{sec:oddsandwiches}, we show that $\mathcal{D}_n^I(x)$ is indeed supported on a set of some special chessboard elements defined in terms of odd sandwiches.
This in turn allows us to factor $\mathcal{D}_n^J(x)$ as the product of $\mathcal{D}_n^{[0,a_0-1]\cup[a_0+1,n-1]}(x)$ and a sign-twisted generating function in type $A$.

\subsection{Chessboard elements}\label{subsec:Chessboardelements}
Chessboard elements were introduced in \cite{KV09} for the symmetric group and generalized in \cite{SV13} to the hyperoctahedral group.
A crucial property of these elements is that the sign-twisted generating function in type $B$ is the same when the support set is restricted to the set of chessboard elements \cite[Lemma 8]{SV13}.
In this subsection, we prove an analogous result for $\mathcal{D}_n^{I}(x)$, so that $\mathcal{D}_n^{I}(x)$ is supported on a relatively small subset of $\mathcal{D}_n^{I}$.

\begin{defn}\label{defn:chessboardelem}
An even signed permutation $\sigma$ in $\mathcal{D}_n$ is called a \emph{chessboard element} if
$i+\sigma(i)\equiv j+\sigma(j)(\bmod\, 2)$ for all $i,j\in[n]$, or equivalently, $\sigma(i)\not\equiv \sigma(i+1)(\bmod\, 2)$ for all $i\in[n-1]$.
\end{defn}

Denote $C_n:=C_{n,0}\cup C_{n,1}$, where
$
C_{n,j}:=\{\sigma\in \mathcal{D}_n\,|\,i+\sigma(i)\equiv j (\bmod\ 2)\ \text{for all}\ i\in[n]\}
$
for $j\in\{0,1\}$.
The set $C_n$ forms a subgroup of $\mathcal{D}_n$, and $C_{n,0}$ forms a subgroup of $C_{n}$.
Moreover, if $n\equiv 1(\bmod\, 2)$, then $C_{n}=C_{n,0}$. For any subset $X\subseteq \mathcal{D}_n$, we write $X^I:= X\cap \mathcal{D}_n^I$.
In particular, we have $C_n^I= C_n\cap\mathcal{D}_n^I$.

Chessboard elements can be used to control cancellation in the sign-twisted generating function $\mathcal{D}_n^I(x)$.

\begin{lemma}\label{lem:DnI=CnIeq}
For any $I\subseteq[0,n-1]$, we have $\mathcal{D}_n^I(x)=\sum_{\sigma\in C_n^I}(-1)^{\ell(\sigma)}x^{L(\sigma)}$.
\end{lemma}
\begin{proof}
Let $\sigma\in \mathcal{D}_n^I\setminus C\mathcal{}_n^I$. Then there exists $i\in[n-1]$ such that
$\sigma^{-1}(i)\equiv\sigma^{-1}(i+1) (\bmod\ 2).$
Let $i$ be chosen to be minimal with this property and write $\sigma'=s_i^D\sigma$.
Then $\sigma'^{-1}(i)=\sigma^{-1}(i+1)$ and  $\sigma'^{-1}(i+1)=\sigma^{-1}(i)$, so that $\sigma'^{-1}(i)\equiv\sigma'^{-1}(i+1) (\bmod\ 2)$,
yielding that $\sigma'\in \mathcal{D}_n^I\setminus C_n^I$.
Thus, the map $\sigma\mapsto \sigma'$ is a fixed-point free involution on $\mathcal{D}_n^I\setminus C_n^I$.
It is easy to see that $\ell(\sigma')=\ell(\sigma)\pm 1$ and $L(\sigma')=L(\sigma)$, so that $(-1)^{\ell(\sigma)}x^{L(\sigma)}+(-1)^{\ell(\sigma')}x^{L(\sigma')}=0$. This implies the assertion.
\end{proof}

A completely analogous argument shows that the following result holds for $\mathcal{S}_n$.
\begin{lemma}\label{lem:SnI=CnI}
For any $I\subseteq[n-1]$, we have $\mathcal{S}_n^I(x)=\sum_{\sigma\in \mathcal{S}_n^I\cap C_n}(-1)^{\ell(\sigma)}x^{L(\sigma)}$.
\end{lemma}
\delete{\begin{proof}
Let $\sigma\in \mathcal{S}_n^I\setminus(\mathcal{S}_n^I\cap C_n)$. Then there exists $i\in[n-1]$ such that
$\sigma^{-1}(i)\equiv\sigma^{-1}(i+1)(\bmod\ 2)$. Let $i$ be minimal with this property and set $\sigma'=s_i\sigma$. Then $\Des(\sigma')=\Des(\sigma)$ since $|\sigma^{-1}(i)-\sigma^{-1}(i+1)|\geq 2$. Hence $\sigma'\in \mathcal{S}_n^I\setminus(\mathcal{S}_n^I\cap C_n)$. Notice that $(\sigma')'=\sigma$ and $\sigma'\neq\sigma$. Moreover, $L(\sigma')=L(\sigma)$ but $\ell(\sigma')=\ell(\sigma)\pm 1$. Thus, for each $\sigma\in \mathcal{S}_n^I\setminus(\mathcal{S}_n^I\cap C_n)$ there exists a unique $\sigma'\in \mathcal{S}_n^I\setminus(\mathcal{S}_n^I\cap C_n)$ such that $(-1)^{\ell(\sigma)}x^{L(\sigma)}+(-1)^{\ell(\sigma')}x^{L(\sigma')}=0$. This implies the assertion.
\end{proof}
}

Let $I=[0, a_0-1]\cup[a_1, b_1]\cup\cdots\cup[a_s,b_s]$ be a subset of $[0,n-1]$.
Recall that for each $\sigma\in \mathcal{D}_n$, $\sigma^I$ is obtained from $\sigma$ by first rearranging the values $|\sigma(1)|, \ldots, |\sigma(a_0)|$ in increasing order,
then changing the sign of the smallest element of $|\sigma(1)|, \ldots, |\sigma(a_0)|$ if the number of negative numbers in $\sigma(1), \ldots, \sigma(a_0)$ is odd,
finally rearranging the values $\sigma(a_i), \ldots, \sigma(b_i+1)$ so that they appear in increasing order in the places $a_i,\cdots, b_i+1$ for all $i\in[s]$,
and leaving all other values unchanged.
For example, let $n=7$, $\sigma=\bar32\bar71\bar54\bar6$, and $\tau=\bar3271\bar546$. If  $I=[0,2]\cup[4,5]$, then
\begin{align*}
\sigma^I=237\bar514\bar6\quad\text{and}\quad \tau^I=\bar237\bar5146.
\end{align*}

We say that an element $\sigma\in \mathcal{D}_n$ is \emph{ascending} if $\sigma(1)<\sigma(2)<\cdots<\sigma(n)$.
Then each $\sigma\in \mathcal{D}_n$ can be decomposed as $\sigma=\sigma^{[n-1]}\sigma_{[n-1]}$, where $\sigma^{[n-1]}\in \mathcal{D}_n^{[n-1]}$ is ascending and
$\sigma_{[n-1]}\in (\mathcal{D}_n)_{[n-1]}\cong \mathcal{S}_n$. Moreover, $\ell(\sigma)=\ell(\sigma^{[n-1]})+\ell(\sigma_{[n-1]})$.
%Indeed, $\sigma^{[n-1]}$ is obtained from $\sigma$ by rearranging the numbers $\sigma(1),\sigma(2),\ldots,\sigma(n)$ in increasing order.
The odd length function $L$ of chessboard elements of $\mathcal{D}_n$ also has additive property under this factorization.
\begin{lemma}\label{lem:L(w)=L(w1)+L(w2)}
If $\sigma\in C_n$, then
$L(\sigma)=L(\sigma^{[n-1]})+L(\sigma_{[n-1]}).$
\end{lemma}
\begin{proof}
For brevity, let $\tau=\sigma^{[n-1]}$.
Then $\tau$ is ascending, and hence $L(\tau)=\onsp(\tau)$.
For any $i<j$ in $[n]$, we have
\begin{align*}
\sigma(i)<\sigma(j)\Leftrightarrow ~ \tau(\sigma_{[n-1]}(i))<\tau(\sigma_{[n-1]}(j))
\Leftrightarrow ~ \sigma_{[n-1]}(i)<\sigma_{[n-1]}(j),
\end{align*}
yielding $\oinv(\sigma)=\oinv(\sigma_{[n-1]})=L(\sigma_{[n-1]})$, since $\sigma_{[n-1]}\in (\mathcal{D}_n)_{[n-1]}\cong \mathcal{S}_n$. But
$L(\sigma)=\onsp(\sigma)+\oinv(\sigma)$, so it remains to show that $\onsp(\tau)=\onsp(\sigma)$.

Suppose now that $\sigma(i_1)<\sigma(i_2)<\cdots<\sigma(i_k)<0$ are all the negative entries of $\sigma$ in one-line notation.
\delete{
\begin{align*}
  \Nsp(\sigma)
=\ \{(j,i)\,|\,1\leq j<i,   \sigma(i)<-|\sigma(j)|\}
\uplus\{(i,j)\,|\,i<j\leq n,  \sigma(i)<-|\sigma(j)|\}.
\end{align*}
But $\sigma(i)<-|\sigma(j)|$ is equivalent to $\sigma(i)<0$ and $\sigma(i)<\sigma(j)<\overline{\sigma(i)}$, so that}
It follows from $\sigma\in C_n$ that
\begin{align*}
  \onsp(\sigma)
= \sum_{s=1}^k|\{j\in[n]\,|\,\sigma(i_s)<\sigma(j)<\overline{\sigma(i_s)}, j\not\equiv i_s (\bmod\, 2)\}|=\sum_{s=1}^k\left\lfloor\frac{|\sigma(i_s)|}{2}\right\rfloor.
\end{align*}
On the other hand, by the definition of $\tau$, we see that $\tau(s)=\sigma(i_s)$ for all $s\in[k]$, and they are all the negative entries of $\tau$ in one-line notation.
Thus, for $s\in[k]$, we have
$$
\{|\tau(s+1)|,|\tau(s+2)|,\ldots,|\tau(s+|\tau(s)|-1)|\}=\{1,2,\ldots,|\tau(s)|-1\}.
$$
It follows that
\begin{align*}
\onsp(\tau)&=\sum_{s=1}^k|\{j\,|\,s<j\leq n, \tau(s)+\tau(j)<0, s\not\equiv j (\bmod\ 2)\}|\\
&=\sum_{s=1}^k|\{j\,|\,s+1\leq j\leq s+|\tau(s)|-1,  s\not\equiv j (\bmod\ 2)\}|\\
&=\sum_{s=1}^k\left\lfloor\frac{|\tau(s)|}{2}\right\rfloor,
\end{align*}
and hence $\onsp(\tau)=\onsp(\sigma)$, so the proof follows.
\end{proof}

For an element $\sigma$ not in $C_n$, it is not necessarily true in general that $L(\sigma)=L(\sigma^{[n-1]})+L(\sigma_{[n-1]})$.
For example, $\sigma=\bar 1\bar 324$ is an element of $\mathcal{D}_4\setminus C_4$ with $L(\sigma)=3$. But $\sigma^{[3]}=\bar 3\bar 124$ and $\sigma_{[3]}=2134$, so that $L(\sigma^{[3]})=L(\sigma_{[3]})=1$. Hence $L(\sigma)\neq L(\sigma^{[3]})+L(\sigma_{[3]})$.

\subsection{Odd sandwiches}\label{sec:oddsandwiches}
To begin with, we define odd sandwiches in a final segment of an even signed permutation, by generalizing \cite[Definition 15]{SV13}.
Recall that the sign function of nonzero real numbers is given by $\sgn(x)=1$ if $x>0$, $\sgn(x)=-1$ if $x<0$.

\begin{defn}\label{def:oddsanwichin[k+1,n]}
	Let $\sigma\in \mathcal{D}_n$, $c\in[n-1]$, and let $r,h\in[n-2]$ such that $h\equiv1(\bmod\ 2)$, $r$ and $r+h+1$ are elements of $\{|\sigma(c)|,\ldots,|\sigma(n)|\}$. Then the pair $(r,h)$ is called an \emph{odd sandwich} in $\sigma(c)\cdots\sigma(n)$ if it satisfies one of the following conditions:
	\begin{enumerate}
%\item\label{item1:def:oddsandwichin[k+1,na]} $r$ and $r+h+1$ are elements of $\{|\sigma(c)|,\ldots,|\sigma(n)|\}$,
\item\label{item1:def:oddsandwichin[k+1,n]} $\sgn(\sigma^{-1}(r))=\sgn(\sigma^{-1}(r+h+1))$ and $\sgn(\sigma^{-1}(r))\neq\sgn(\sigma^{-1}(s))$ for all $s\in[r+1,r+h]\cap\{|\sigma(c)|,\ldots,|\sigma(n)|\}$,
\item\label{item2:def:oddsandwichin[k+1,n]}  $r=\min\{|\sigma(c)|,\ldots,|\sigma(n)|\}$, $\sgn(\sigma^{-1}(r))\neq\sgn(\sigma^{-1}(r+h+1))$ and $\sgn(\sigma^{-1}(r))=\sgn(\sigma^{-1}(s))$ for all $s\in[r+1,r+h]\cap\{|\sigma(c)|,\ldots,|\sigma(n)|\}$.
	\end{enumerate}
%We say that $\sigma(c)\cdots\sigma(n)$ \emph{has an odd sandwich} if there exists an odd sandwich in $\sigma(c)\cdots\sigma(n)$.
\end{defn}

For example, if $\sigma=\overline{4}\,\overline{3}52\overline{1}6\overline{7}$, then $(3,3)$ is an odd sandwich in $\overline{3}52\overline{1}6\overline{7}$.
If $\sigma=\overline{1}56\overline{3}\,\overline{2}7\overline{4}$, then $(2,3)$ is an odd sandwich in $6\overline{3}\,\overline{2}7\overline{4}$.
Taking $c=1$, Definition \ref{def:oddsanwichin[k+1,n]} recovers the notion of odd sandwiches in signed permutations (i.e., \cite[Definition 15]{SV13}).
%It is easy to see that for any $\sigma,\tau\in \mathcal{D}_n$, if $\{\sigma(i)\,|\,i\in [c,n]\}=\{\tau(i)\,|\,i\in [c,n]\}$, then $\sigma(c)\cdots\sigma(n)$ and $\tau(c)\cdots\tau(n)$ have %the same set of odd sandwiches.

For any $c\in[n-1]$, we write
	$$
	H_{n,c}:=\{\sigma\in C_n\,|\, \sigma(c)\cdots\sigma(n) \text{ has no odd sandwiches}\}.
	$$
The following result says that $\mathcal{D}_n^I(x)$ is supported on a set of the form $H_{n,c}^I:=H_{n,c}\cap \mathcal{D}_n^I$ for some suitable $c\in[n-1]$.
\begin{lemma}\label{lem:CnI=HnILl}
Let $I\subseteq [0,n-1]$ with $I_0=[0,a_0-1]$. If $a_0\in[2,n-1]$, then
	$$
\mathcal{D}_n^I(x)=\sum_{\sigma\in H_{n,a_0+1}^I}(-1)^{\ell(\sigma)}x^{L(\sigma)}.
	$$
\end{lemma}
\begin{proof}
By Lemma \ref{lem:DnI=CnIeq}, it suffices to show that
\begin{align}\label{eq:cnI=Hna0+1I}
\sum_{\sigma\in C_n^I\setminus H_{n,a_0+1}^I}(-1)^{\ell(\sigma)}x^{L(\sigma)}=0.
\end{align}
Notice that
	$
	C_n^I\setminus H_{n,a_0+1}^I=(\mathcal{D}_n^I\cap C_n)\setminus(\mathcal{D}_n^I\cap H_{n,a_0+1})=\mathcal{D}_n^I\cap(C_n\setminus H_{n,a_0+1}).
	$
Define a map $\varphi: C_n^I\setminus H_{n,a_0+1}^I\rightarrow  C_n^I\setminus H_{n,a_0+1}^I$ by sending $\sigma$ to
\begin{align*}
\sigma'=\begin{cases}
(\barr{r}, \barr{r+h+1})(r, r+h+1)\sigma, &{\rm if}\ (r,h)\ {\rm satisfies\ Definition}\ \ref{def:oddsanwichin[k+1,n]}\eqref{item1:def:oddsandwichin[k+1,n]},\\
(\barr{r}, {r+h+1})(r, \barr{r+h+1})\sigma, &{\rm if}\ (r,h)\ {\rm satisfies\ Definition}\ \ref{def:oddsanwichin[k+1,n]}\eqref{item2:def:oddsandwichin[k+1,n]},
\end{cases}
\end{align*}
where $(r,h)$ is the minimum odd sandwich in $\sigma(a_0+1)\cdots\sigma(n)$ under the lexicographical order.
Then $r\equiv r+h+1(\bmod\,2)$. Since $a_0\in[2,n-1]$, we see that $\sigma\in C_n\backslash H_{n,a_0+1}$ is equivalent to $\sigma'\in C_n\backslash H_{n,a_0+1}$.
It follows from Definition \ref{def:oddsanwichin[k+1,n]} that $\sigma\in \mathcal{D}_n^{I}$ if and only if $\sigma'\in \mathcal{D}_n^{I}$, from which we see that $\sigma\in C_n^I\setminus H_{n,a_0+1}^I$ is equivalent to $\sigma'\in C_n^I\setminus H_{n,a_0+1}^I$.
Thus, $\varphi$ is well defined.
Since the odd sandwich $(r,h)$ is chosen to be lexicographically minimal, the $(r,h)$ is unique for each $\sigma$.
Thus the map $\sigma\mapsto \sigma'$ is an involution on the set $C_n^I\setminus H_{n,a_0+1}^I$ without fixed points.

Since $h$ is odd, we have $\oinv(\sigma)=\oinv(\sigma')$ and $\onsp(\sigma)=\onsp(\sigma')$, so that $L(\sigma')=L(\sigma)$. But $\ell(\sigma')=\ell(\sigma)\pm 1$, so $(-1)^{\ell(\sigma)}x^{L(\sigma)}+(-1)^{\ell(\sigma')}x^{L(\sigma')}=0$. Thus, \eqref{eq:cnI=Hna0+1I} holds,
and the proof follows.
\end{proof}

The following result reveals a close relationship between chessboard elements and odd sandwiches.
\begin{lemma}\label{lem:relationoscs}
	Let $\sigma\in \mathcal{D}_n$ and $c\in[n-1]$ such that $\sigma(c)<\cdots<\sigma(n)$. Then $\sigma(c)\cdots\sigma(n)$ has an odd sandwich if and only if there exists $i\in[c,n-1]$ such that
	$\sigma(i)\equiv \sigma(i+1)(\bmod\, 2)$.
	In particular, if $\sigma$ is ascending, then $\sigma$ has no odd sandwiches if and only if $\sigma$ is a chessboard element.
\end{lemma}
\begin{proof}
	Assume that there exists $i\in[c,n-1]$ such that $\sigma(i)\equiv \sigma(i+1)(\bmod\ 2)$.
	Let $h=||\sigma(i)|-|\sigma(i+1)||-1$. It is clear that $h\in[n-2]$ is odd.  Let $r=\min\{|\sigma(i)|, |\sigma(i+1)|\}$.
	Then $r\in[n-2]$ and $r+h+1=\max\{|\sigma(i)|, |\sigma(i+1)|\}.$
	If $\sgn(\sigma^{-1}(r))=\sgn(\sigma^{-1}(r+h+1))$, then $\sgn(\sigma^{-1}(r))\neq\sgn(\sigma^{-1}(s))$ for all $s\in[r+1,r+h]\cap\{|\sigma(c)|,\ldots,|\sigma(n)|\}$, since $\sigma(c)<\cdots<\sigma(n)$. Thus, $(r,h)$ is an odd sandwich in $\sigma(c)\cdots\sigma(n)$.
	If $\sgn(\sigma^{-1}(r))\neq\sgn(\sigma^{-1}(r+h+1))$, then $\sigma(i)<0<\sigma(i+1)$, so that $r={\rm min}\{|\sigma(c)|,\ldots,|\sigma(n)|\}$. Notice that $\sigma(c)<\cdots<\sigma(n)$, so
	$\sgn(\sigma^{-1}(r))=\sgn(\sigma^{-1}(s))$ for all $s\in[r+1,r+h]\cap\{|\sigma(c)|,\ldots,|\sigma(n)|\}$, and hence $(r,h)$ is an odd sandwich in $\sigma(c)\cdots\sigma(n)$.
	Therefore, in this case $\sigma(c)\cdots\sigma(n)$ has an odd sandwich.
	
	For the other direction, let $(r,h)$ be an odd sandwich of $\sigma(c)\cdots\sigma(n)$, and let $p=\sigma^{-1}(r)$ and $q=\sigma^{-1}(r+h+1)$. Then $|p|$ and $|q|$ are elements of $[c,n]$.
First consider the case $p>0$. We have $p\in[c,n]$ and $\sigma(p)=r$.
If $(r,h)$ satisfies Definition \ref{def:oddsanwichin[k+1,n]}\eqref{item1:def:oddsandwichin[k+1,n]}, then $q>0$ and $\sigma(q)=r+h+1$,
but $\sigma^{-1}(s)<0$ for all $s\in[r+1,r+h]\cap\{|\sigma(c)|,\ldots,|\sigma(n)|\}$.
Since $\sigma(c)<\cdots<\sigma(n)$, it follows that $q=p+1$. If  $(r,h)$ satisfies Definition \ref{def:oddsanwichin[k+1,n]}\eqref{item2:def:oddsandwichin[k+1,n]}, then $q<0$ and $\sigma(q)=r+h+1$,
but $\sigma^{-1}(s)>0$ for all $s\in[r+1,r+h]\cap\{|\sigma(c)|,\ldots,|\sigma(n)|\}$.
So we deduce from $\sigma(c)<\cdots<\sigma(n)$ that $q=p-1$.
Similarly, if $p<0$, then a completely analogous argument shows that $q=p\pm 1$. Therefore, $|p-q|=1$, that is,
	$||\sigma^{-1}(r)|-|\sigma^{-1}(r+h+1)||=1$.
	Hence, $i=\min\{|\sigma^{-1}(r)|, |\sigma^{-1}(r+h+1)|\}$ is the desired integer in $[c,n-1]$ such that $\sigma(i)\equiv \sigma(i+1)(\bmod\, 2)$.
\end{proof}

Given compressed subsets $I=[0,a_0-1]\cup[a_0+1,a_1-1]\cup\cdots\cup[a_{s-1}+1,n-1]$ and $J=[0,a_0-1]\cup[a_0+1,n-1]$, where $a_0\geq2$, we
now give a decomposition of the set $H_{n,a_0+1}^I$ by using the parabolic factorizations of the form
$\sigma=\sigma^{J}\sigma_{J}$, from which we deduce an elegant product formula for $\mathcal{D}_n^I(x)$ in Proposition \ref{lem:DnI=DnJxSn-aKeq}.

Given positive integers $p$ and $q$, let $\mathcal{D}_p\times\mathcal{D}_q$ denote the subgroup of $\mathcal{D}_{p+q}$, where $\mathcal{D}_{p}$ permutes $[p]_{\pm}$ and
$\mathcal{D}_{q}$ permutes $[\overline{p+q},\overline{p+1}]\cup[p+1,p+q]$.
So for any $\sigma\in \mathcal{D}_p$ and $\tau\in \mathcal{D}_q$, we have
$$
\sigma\times\tau=(\sigma(1),\ldots,\sigma(p),\tau(1)+\sgn(\tau(1))p,\ldots,\tau(q)+\sgn(\tau(q))p).
$$
For example, $12\times 34\bar2\bar1=1256\bar4\bar3$.

\begin{lemma}\label{lem:HnI=(HnJ)(SnKcupCn)1}
Let $I=[0,a_0-1]\cup[a_0+1,a_1-1]\cup\cdots\cup[a_{s-1}+1,n-1]$ be compressed with $a_0\in[2,n-1]$,
and let $J=[0,a_0-1]\cup[a_0+1,n-1]$ and $K=[1,a_1-a_0-1]\cup\cdots\cup[a_{s-1}-a_0+1,n-a_0-1]$. Then
\begin{align*}
H_{n,a_0+1}^J=\{\sigma^{J}\,|\, \sigma\in H_{n,a_0+1}^I\},\quad (\{1_{a_0}\}\times \mathcal{S}_{n-a_0}^K)\cap C_n=\{\sigma_{J}\,|\, \sigma\in H_{n,a_0+1}^I\},
\end{align*}
and $H_{n,a_0+1}^I=H_{n,a_0+1}^J\cdot ((\{1_{a_0}\}\times \mathcal{S}_{n-a_0}^K)\cap C_n)$.
Here $1_{a_0}$ is the identity element of $\mathcal{D}_{a_0}$.
\end{lemma}

\begin{proof}
%We assume that $a_0\in[2,n-2]$. The proof for the case $a_0=0$ is similar and will be omitted.
Let $\sigma\in H_{n,a_0+1}^J$ and $\tau\in (\{1_{a_0}\}\times \mathcal{S}_{n-a_0}^K)\cap C_n$.
Since $C_n$ is a subgroup of $\mathcal{D}_n$, we see that $\sigma\tau\in C_n$.
%Then $\sigma(a_0+1)<\cdots<\sigma(n)$ and $\{\sigma\tau(i)\,|\,i\in[a_0+1,n]\}=\{\sigma(i)\,|\,i\in[a_0+1,n]\}$, so that $\sigma\tau\in H_{n,a_0+1}$.
Notice that $\{\tau(i)\,|\,i\in[a_0+1,n]\}=[a_0+1,n]$. So $\sigma\tau\in H_{n,a_0+1}$.
Since $\sigma\in \mathcal{D}_n^J$ and $[a_0+1,n-1]\subseteq J$, we have $\sigma(a_0+1)<\cdots<\sigma(n)$. But $a_0\in[2,n-1]$, we see that $\sigma\tau\in \mathcal{D}_n^I$, from which we get $\sigma\tau\in H_{n,a_0+1}^I$.
Hence, $H_{n,a_0+1}^J\cdot ((\{1_{a_0}\}\times \mathcal{S}_{n-a_0}^K)\cap C_n)\subseteq H_{n,a_0+1}^I$.
It is easy to see that $\sigma=(\sigma\tau)^J$ and $\tau=(\sigma\tau)_J$. Thus,
\begin{align*}
H_{n,a_0+1}^J\subseteq\{\sigma^{J}\,|\, \sigma\in H_{n,a_0+1}^I\}\quad \text{and}\quad (\{1_{a_0}\}\times \mathcal{S}_{n-a_0}^K)\cap C_n\subseteq\{\sigma_{J}\,|\, \sigma\in H_{n,a_0+1}^I\}.
\end{align*}

Let $\sigma\in H_{n,a_0+1}^I$. It is enough to show that $\sigma^J\in H_{n,a_0+1}^J$ and $\sigma_J\in (\{1_{a_0}\}\times \mathcal{S}_{n-a_0}^K)\cap C_n$.
Clearly, $\sigma^J\in \mathcal{D}_n^J$ and $\sigma_J\in \{1_{a_0}\}\times \mathcal{S}_{n-a_0}^K$.
Moreover, $\sigma^J(a_0+1)\cdots\sigma^J(n)$ is obtained from $\sigma$ by rearranging $\sigma(a_0+1),\ldots,\sigma(n)$
in increasing order. Thus, $\sigma(a_0+1)\cdots\sigma(n)$ and $\sigma^J(a_0+1)\cdots\sigma^J(n)$ have the same set of odd sandwiches, so that $\sigma^J(a_0+1)\cdots\sigma^J(n)$ has no
odd sandwiches.
It remains to show that $\sigma^J$ and $\sigma_J$ are elements of $C_n$.

It follows from Lemma \ref{lem:relationoscs} that $\sigma^J(i)\not\equiv \sigma^J(i+1)(\bmod\, 2)$ for all $i\in[a_0+1,n-1]$.
Now $\sigma\in H_{n,a_0+1}^I$ implies $\sigma\in C_n$, and hence $\sigma(i)\not\equiv \sigma(i+1)(\bmod\,2)$ for all $i\in[n-1]$.
But $\sigma^J(i)=\sigma(i)$ for all $i\in[a_0]$, so  $\sigma^J(i)\not\equiv \sigma^J(i+1)(\bmod\, 2)$ for all $i\in[a_0-1]$.
Notice that
\begin{align*}
\sigma^J(a_0+1)={\rm min}\{\sigma(i)\,|\,i\in[a_0+1,n]\}={\rm min}\{\sigma(a_i+1)\,|\,i\in[0,s-1]\}.
\end{align*}
Since $I$ is compressed, we have $a_0\equiv a_1\equiv\cdots\equiv a_{s-1}\equiv n(\bmod\, 2)$.
Thus, there exists $i\in[0,s-1]$ such that
$$\sigma^J(a_0+1)= \sigma(a_i+1)\equiv\sigma(a_0+1)\not\equiv\sigma(a_0)=\sigma^J(a_0)(\bmod\,2).$$
Hence, $\sigma^J(i)\not\equiv\sigma^J(i+1)$ for all $i\in[n-1]$, that is, $\sigma^J\in C_n$.
Since $C_n$ is a subgroup of $\mathcal{D}_n$, we see that $\sigma_J=(\sigma^J)^{-1}\sigma\in C_n$,
completing the proof.
\end{proof}

\begin{lemma}\label{lem:H_n^Ifactor_L=L+L}
	Let  $I=[0,a_0-1]\cup[a_0+1,a_1-1]\cup\cdots\cup[a_{s-1}+1,n-1]$ be compressed with $a_0\in[2,n-1]$. Then  $L(\sigma)=L(\sigma^J)+L(\sigma_J)$ for any  $\sigma\in H_{n,a_0+1}^I$, where
	$J=[0,a_0-1]\cup[a_0+1,n-1]$.
\end{lemma}
\begin{proof}
	Let $\sigma\in H_{n,a_0+1}^I$. By Lemma \ref{lem:HnI=(HnJ)(SnKcupCn)1}, $\sigma^J\in H_{n,a_0+1}^J$ and $\sigma_J\in(\{1_{a_0}\}\times \mathcal{S}_{n-a_0}^K)\cap C_n$, where $K=[1,a_1-a_0-1]\cup\cdots\cup[a_{s-1}-a_0+1,n-a_0-1]$.
	Thus,	
	\begin{align*}
		L(\sigma) =&\ |\{(i,j)\in[a_0]\times[a_0+1,n]\,|\,\sigma(i)>\sigma(j),i\not\equiv j(\bmod\ 2)\}|\\
		&\ +|\{(i,j)\in[a_0+1,n]^2\,|\,i<j,\sigma(i)>\sigma(j),i\not\equiv j(\bmod\ 2)\}|\\
		&\ +|\{(i,j)\in[a_0]\times[a_0+1,n]\,|\,\sigma(i)+\sigma(j)<0,i\not\equiv j(\bmod\ 2)\}|\\
		&\ +|\{(i,j)\in[a_0]^2\cup [a_0+1,n]^2\,|\,i<j,\sigma(i)+\sigma(j)<0,i\not\equiv j(\bmod\ 2)\}|,\\
		L(\sigma^J) =&\ |\{(i,j)\in[a_0]\times[a_0+1,n]\,|\,\sigma^J(i)>\sigma^J(j),i\not\equiv j(\bmod\ 2)\}|\\
		&\ +|\{(i,j)\in[a_0]\times[a_0+1,n]\,|\,\sigma^J(i)+\sigma^J(j)<0,i\not\equiv j(\bmod\ 2)\}|\\
		&\ +|\{(i,j)\in[a_0]^2\cup [a_0+1,n]^2\,|\,i<j,\sigma^J(i)+\sigma^J(j)<0,i\not\equiv j(\bmod\ 2)\}|,\\
		L(\sigma_J) =&\ |\{(i,j)\in[a_0+1,n]^2\,|\,i<j,\sigma_J(i)>\sigma_J(j),i\not\equiv j(\bmod\ 2)\}|.
	\end{align*}
	Since $\sigma(i)>\sigma(j)$ if and only if $\sigma_J(i)>\sigma_J(j)$ for any $(i,j)\in[a_0+1,n]^2$,
the second term on the right-hand side of $L(\sigma)$ is equal to $L(\sigma_J)$.
	Notice that $$\sigma^J(n)={\rm max} \{\sigma(i)\,|\, i\in[a_0+1,n]\}={\rm max} \{\sigma(a_i)\,|\, i\in[s]\},$$
	where $a_s=n$.
	Assume that $\sigma^J(n)=\sigma(a_r)$ for some $r\in[s]$.
	By the definition of compressed subsets, we have $a_r\equiv n(\bmod\ 2)$.
	Since $\sigma,\sigma^J\in C_n$, there follows
	\begin{align*}
		\sigma^J(i)+i\equiv \sigma^J(n)+n\equiv \sigma(a_r)+a_r(\bmod\ 2)
	\end{align*}
	for all $i\in[a_0+1,n]$. In particular, for any $i\in[a_0+1,n]$, we have $(\sigma^J)^{-1}(\sigma(i))\in [a_0+1,n]$,
	so that
	\begin{align*}
		(\sigma^J)^{-1}(\sigma(i))+\sigma(i)\equiv \sigma(a_r)+a_r\equiv \sigma(i)+i(\bmod\ 2),
	\end{align*}
	where the second congruence follows from $\sigma\in C_n$. Hence, $(\sigma^J)^{-1}(\sigma(i))\equiv i(\bmod\ 2)$ for all $i\in[a_0+1,n]$.
	But $\sigma^J(i)=\sigma(i)$ for all $i\in[a_0]$, so $(\sigma^J)^{-1}(\sigma(i))\equiv i(\bmod\ 2)$ for all $i\in[n]$. Clearly,
	$\sigma(i)=\sigma^J((\sigma^J)^{-1}(\sigma(i)))$ for all $i\in[n]$,
	so the sum of the first, third and forth terms on the right-hand side of $L(\sigma)$ is equal to $L(\sigma^J)$.
	Hence, $L(\sigma)=L(\sigma^J)+L(\sigma_J)$.
\end{proof}

\begin{prop}\label{lem:DnI=DnJxSn-aKeq}
Let $I=[0,a_0-1]\cup[a_0+1,a_1-1]\cup\cdots\cup[a_{s-1}+1,n-1]$ be a compressed subset of $[0,n-1]$, and let $J=[0,a_0-1]\cup[a_0+1,n-1]$. If $a_0\in[2,n-1]$, then
	\begin{align*}
		\mathcal{D}_n^I(x)=\left[\begin{matrix}\frac{n-a_0}{2}\\ \frac{a_1-a_0}{2},\ldots,\frac{n-a_{s-1}}{2}\end{matrix}\right]_{x^2}
\mathcal{D}_n^J(x).
	\end{align*}
\end{prop}
\delete{\begin{proof}
By Lemma \ref{lem:CnI=HnILl}, $\mathcal{D}_n^I(x)=\sum\limits_{\sigma\in H_{n,a_0+1}^I}(-1)^{\ell(\sigma)}x^{L(\sigma)}$.
Let $K=[1,a_1-a_0-1]\cup\cdots\cup[a_{s-1}-a_0+1,a_s-a_0-1]$. It follows from  Lemmas \ref{lem:L(w)=L(w1)+L(w2)}, \ref{lem:HnI=(HnJ)(SnKcupCn)1}
and Lemmas \ref{lem:SnI=CnI}, \ref{lem:CnI=HnILl} that
	\begin{align*}
\mathcal{D}_n^I(x)=\Big(\sum_{\sigma\in H_{n,a_0+1}^J}(-1)^{\ell(\sigma)}x^{L(\sigma)}\Big)\Big(\sum_{\sigma\in (\{1_{a_0}\}\times \mathcal{S}_{n-a_0}^K)\cap C_n}(-1)^{\ell(\sigma)}x^{L(\sigma)}\Big)
		=\mathcal{D}_n^J(x) \mathcal{S}_{n-a_0}^K(x),
	\end{align*}
%and by  Lemmas \ref{lem:SnI=CnI} and \ref{lem:CnI=HnILl}, $\mathcal{D}_n^I(x)=\mathcal{D}_n^J(x) \mathcal{S}_{n-a_0}^K(x)$,
so the proof follows from Theorem \ref{theo:SnI(x)}.
\end{proof}
}
\begin{proof}
Let $K=[1,a_1-a_0-1]\cup\cdots\cup[a_{s-1}-a_0+1,n-a_0-1]$. Then
	\begin{align*}
\mathcal{D}_n^I(x)=&\sum_{\sigma\in H_{n,a_0+1}^I}(-1)^{\ell(\sigma)}x^{L(\sigma)}&&{\rm by\ Lemma\ \ref{lem:CnI=HnILl}}\\
		=&\, \Big(\sum_{\sigma\in H_{n,a_0+1}^J}(-1)^{\ell(\sigma)}x^{L(\sigma)}\Big)\Big(\sum_{\sigma\in (\{1_{a_0}\}\times \mathcal{S}_{n-a_0}^K)\cap C_n}(-1)^{\ell(\sigma)}x^{L(\sigma)}\Big)&&
{\rm by\ Lemmas\ \ref{lem:HnI=(HnJ)(SnKcupCn)1}\ and\ \ref{lem:H_n^Ifactor_L=L+L}}\\
		=& \,\mathcal{D}_n^J(x) \mathcal{S}_{n-a_0}^K(x), &&{\rm by\ Lemmas\ \ref{lem:SnI=CnI}\ and\ \ref{lem:CnI=HnILl}}
	\end{align*}
so the proof follows from Theorem \ref{theo:SnI(x)}.
\end{proof}

\section{Closed product formulas for $\mathcal{D}_n^{I}(x)$}\label{sec:formudnixfull}

In this section, we give closed product formulas of $\mathcal{D}_n^{I}(x)$ for all $I\subseteq[0,n-1]$, and thereby verify Conjectures \ref{conj:DnI0i} and \ref{conj:DnI(x)}.
We first consider the case where $I=[0,a_0-1]\cup[a_0+1,n-1]$ is compressed with $|I_0|\neq1$ in the first three subsections. %in Subsections \ref{I0=0emp}, \ref{I0=oddpositiveint} and \ref{I0=evenpositiveint}.
There are three cases that need to be discussed individually: $|I_0|=0$, $|I_0|$ is an odd positive integer and $|I_0|$ is an even positive integer.
In particular, for the case $|I_0|=0$, we establish directly
formulas of $\mathcal{D}_n^{I}(x)$ for all subsets (not necessary compressed) $I$ of $[n-1]$.
We then establish our main theorem (i.e., Theorem \ref{theo:finaltheorem}) in Subsection \ref{mainthm}.

\subsection{When $|I_0|=0$}\label{I0=0emp}

 We now give a closed product formula for $\mathcal{D}_n^{[n-1]}(x)$.
\begin{lemma}\label{prop:Dn[n-1](x)}
For any positive integer $n$, we have
$$
\mathcal{D}_n^{[n-1]}(x)=\prod_{i=2}^n\left(1+(-1)^{i-1}x^{\left\lfloor\frac{i}{2}\right\rfloor}\right).
$$
\end{lemma}
\begin{proof}
The proof is by induction on $n$. Since
$\mathcal{D}_1^{[0]}(x)=1,$
the assertion is true for $n=1$.
Now assume that the desired identity is true for $n-1$, where $n\geq 2$. Since any $\sigma\in \mathcal{D}_n^{[n-1]}$ is ascending, we have the following disjoint union
\begin{align*}
\mathcal{D}_n^{[n-1]}=\{\sigma\in \mathcal{D}_n^{[n-1]}\,|\,\sigma(1)=\barr{n}\}\cup\{\sigma\in \mathcal{D}_n^{[n-1]}\,|\,\sigma(n)=n\}.
\end{align*}
For any $\sigma\in \mathcal{D}_n^{[n-1]}$, if we let $\sigma'=\sigma(2)\cdots\sigma(n)$ and $\sigma''=(1,\overline{1})\sigma'$,
then the map $\sigma\mapsto \sigma''$ is a bijection from $\{\sigma\in \mathcal{D}_n^{[n-1]}\,|\,\sigma(1)=\barr{n}\}$ to $\mathcal{D}_{n-1}^{[n-2]}$.
Moreover, $\ell(\sigma)=\ell(\sigma'')+n-1$ and $L(\sigma)=L(\sigma'')+\left\lfloor\frac{n}{2}\right\rfloor$.
Thus,
\begin{align*}
\mathcal{D}_n^{[n-1]}(x)=&\ \sum_{{\sigma\in \mathcal{D}_n^{[n-1]},\, \sigma(n)=n}}(-1)^{\ell(\sigma)}x^{L(\sigma)}
+\sum_{{\sigma\in \mathcal{D}_n^{[n-1]},\, \sigma(1)=\bar n}}(-1)^{\ell(\sigma)}x^{L(\sigma)}\\
=&\ \mathcal{D}_{n-1}^{[n-2]}(x)+(-1)^{n-1}x^{\left\lfloor\frac{n}{2}\right\rfloor}\sum_{\sigma\in  \mathcal{D}_{n-1}^{[n-2]}}(-1)^{\ell(\sigma)}x^{L(\sigma)}\\
=&\ \left(1+(-1)^{n-1}x^{\left\lfloor\frac{n}{2}\right\rfloor}\right) \mathcal{D}_{n-1}^{[n-2]}(x),
\end{align*}
and the proof follows by induction.
\end{proof}

When $I$ is a subset of $[n-1]$, the following result gives a parabolic decomposition to the set $C_n^I$,
from which we obtain a formula of $\mathcal{D}_n^I(x)$. Recall the definition of $C_I$ from Equation \eqref{eq:m_I}.

\begin{prop}\label{prop:DnI(x)-Iin[n-1]}
	Let $I\subseteq [n-1]$. Then
	%$\mathcal{D}_n^I(x)= \mathcal{D}_n^{[n-1]}(x) \mathcal{S}_n^{I}(x)$, so that
	\begin{align*}%\label{eq:DnI(x)[n-1]}
		\mathcal{D}_n^I(x)
		=C_I\prod_{i=2}^n\left(1+(-1)^{i-1}x^{\left\lfloor\frac{i}{2}\right\rfloor}\right)\prod_{i=2m_I+2}^n\left(1+(-1)^{i-1}x^{\left\lfloor\frac{i}{2}\right\rfloor}\right).
	\end{align*}
\end{prop}
\begin{proof}
	By Lemma \ref{lem:DnI=CnIeq}, $\mathcal{D}_n^I(x)= \sum_{\sigma\in C_n^I}(-1)^{\ell(\sigma)}x^{L(\sigma)}$. Thus, it suffices to show that
	$C_n^I=C_n^{[n-1]}\cdot(\mathcal{S}_n^I\cap C_n)$, where
	\begin{align*}
		C_n^{[n-1]}=\{\sigma^{[n-1]}\,|\, \sigma\in C_n^I\}\quad\text{and}\quad \mathcal{S}_n^I\cap C_n=\{\sigma_{[n-1]}\,|\, \sigma\in C_n^I\},
	\end{align*}
	from which we obtain
	\begin{align*}
		\mathcal{D}_n^I(x)&=\Big(\sum_{\sigma\in C_n^{[n-1]}}(-1)^{\ell(\sigma)}x^{L(\sigma)}\Big)\Big(\sum_{\sigma\in \mathcal{S}_n^I\cap C_n}(-1)^{\ell(\sigma)}x^{L(\sigma)}\Big)
		&&{\rm by\ Lemma\  \ref{lem:L(w)=L(w1)+L(w2)} }\\
		&= \mathcal{D}_n^{[n-1]}(x) \mathcal{S}_n^{I}(x).&&{\rm by\ Lemmas\ \ref{lem:DnI=CnIeq}\ and\  \ref{lem:SnI=CnI}}
	\end{align*}
		Hence, the desired identity  follows from Theorem \ref{theo:SnI(x)} and Lemma \ref{prop:Dn[n-1](x)}.
	
	Let $\delta\in C_n^{[n-1]}$ and $\tau\in \mathcal{S}_n^I\cap C_n$.
	Then $\delta(1)<\cdots<\delta(n)$ and $\delta\tau\in \mathcal{D}_n^I\cap C_n=C_n^I$, which yields that $C_n^{[n-1]}\cdot(\mathcal{S}_n^I\cap C_n)\subseteq C_n^I$. But $\delta=(\delta\tau)^{[n-1]}$ and $\tau=(\delta\tau)_{[n-1]}$, so
	\begin{align*}
		C_n^{[n-1]}\subseteq\{\sigma^{[n-1]}\,|\, \sigma\in C_n^I\}\quad\text{and}\quad \mathcal{S}_n^I\cap C_n\subseteq\{\sigma_{[n-1]}\,|\, \sigma\in C_n^I\}.
	\end{align*}
	Let $\sigma\in C_n^I$. Then  $\sigma^{[n-1]}(1)<\cdots<\sigma^{[n-1]}(n)$, so that $\sigma_{[n-1]}\in\mathcal{S}_n^I$ and that $\sigma^{[n-1]}$ has no odd sandwiches.
	By Lemma \ref{lem:relationoscs}, $\sigma^{[n-1]}\in C_n\cap\mathcal{D}_n^{[n-1]}=C_n^{[n-1]}$,
	so that $\sigma_{[n-1]}=(\sigma^{[n-1]})^{-1}\sigma\in C_n$, since $C_n$ is a group.
	Thus,
	\begin{align*}
		C_n^{[n-1]}\supseteq\{\sigma^{[n-1]}\,|\, \sigma\in C_n^I\}\quad\text{and}\quad \mathcal{S}_n^I\cap C_n\supseteq\{\sigma_{[n-1]}\,|\, \sigma\in C_n^I\},
	\end{align*}
	and hence $C_n^I\subseteq C_n^{[n-1]}\cdot(\mathcal{S}_n^I\cap C_n)$,
	completing the proof.
\end{proof}

\subsection{When $|I_0|$ is an odd positive integer}\label{I0=oddpositiveint}
In this subsection, we give an explicit product formula of $\mathcal{D}_n^I(x)$ when $I$ is a subset of $[0,n-1]$ and $|I_0|\geq2$ is an odd positive integer.
First, for each even signed permutation we introduce the notion of $k$-odd sandwiches, which is a natural analogue of the concept of odd sandwiches.
Then we restrict the supporting set $C_n^I$ of $\mathcal{D}_n^I(x)$ to a relatively small set defined in terms of $k$-odd sandwiches.
By using a parabolic factorization of the form
$\sigma=\sigma^{[0,|I_0|-1]}\sigma_{[0,|I_0|-1]}$, we show that the function $L$ on the smaller supporting set
is additive with respect to this kind of factorization, from which we finally deduce a formula of $\mathcal{D}_n^I(x)$ in Proposition \ref{prop:noddIcompressed}.

\begin{defn}\label{def:k-oddsandwich}
	Let $\sigma\in \mathcal{D}_n$, and let $r,h\in[n-2]$, where $h\equiv1(\bmod\,2)$. For a certain $k\in[n-1]$, the pair $(r,h)$ is called a \emph{$k$-odd sandwich} in $\sigma$ if it satisfies $|\sigma^{-1}(r)|\leq k$, $|\sigma^{-1}(r+h+1)|\leq k$, and $|\sigma^{-1}(r+i)|>k$ for all $i\in[h]$.
\end{defn}

For example, $(3,3)$ is a $4$-odd sandwich in $23\overline{1}74\overline{6}5$.
Given $k\in [n-1]$, let $$T_n(k)=\{\sigma\in C_n\,\big|\,\sigma \text{ has no $k$-odd sandwiches}\}.$$
Since we only consider the set $T_n(a_0)$ for some set $I\subset[0,n-1]$ with $|I_0|=a_0$, for notational simplicity we write $T_n$ for $T_n(a_0)$.

\begin{lemma}\label{lem:eqDnI=TnI}
Let $I=[0,a_0-1]\cup[a_0+1,n-1]$ with $a_0\in[2,n-2]$. Then for any $J\in\{I, I\setminus\{0\}\}$, we have
$
\mathcal{D}_n^J(x)=\sum_{\sigma\in T_n^{J}}(-1)^{\ell(\sigma)}x^{L(\sigma)}.
$
\end{lemma}
\begin{proof}
By Lemma \ref{lem:DnI=CnIeq}, it suffices to show that
\begin{align}\label{eq:TnI=CnI}
\sum_{\sigma\in C_n^{J}\setminus T_n^{J}}(-1)^{\ell(\sigma)}x^{L(\sigma)}=0.
\end{align}
Let $\sigma\in C_n^{J}\setminus T_n^{J}$. Then $\sigma$ has $a_0$-odd sandwiches. Let $(r,h)$ be the minimum $a_0$-odd sandwich in the lexicographical order.
Since $\sigma\in C_n$ and $h\equiv1(\bmod\,2)$, it follows that $|\sigma^{-1}(r)|\equiv|\sigma^{-1}(r+h+1)|(\bmod\,2)$, and hence there exists
$k\in[a_0]$ such that
$$\min\{\sigma(|\sigma^{-1}(r)|),\sigma(|\sigma^{-1}(r+h+1)|)\}< \sigma(k)<\max\{\sigma(|\sigma^{-1}(r)|),\sigma(|\sigma^{-1}(r+h+1)|)\}.$$
But $|\sigma^{-1}(r+i)|>a_0$ for all $i\in[h]$, so $\sigma(k)\not\in[\overline{r+h},\overline{r+1}]\cup[r+1,r+h]$, and hence $\sigma^{-1}(r)$ and $\sigma^{-1}(r+h+1)$ have different signs.

If $J=I$, then we must have $C_n^J=T_n^J$; otherwise, since $|\sigma(1)|<\sigma(2)<\cdots<\sigma(a_0)$ and $|\sigma^{-1}(r+i)|>a_0$ for all $i\in[h]$, it follows that $\sigma(1)=\overline{r}$ and $\sigma(2)=r+h+1$.
Hence  $\sigma(1)\equiv \sigma(2)(\bmod\,2)$,
contradicting $\sigma\in C_n$. Thus, $C_n^J=T_n^J$, that is, $C_n^J\backslash T_n^J=\emptyset$, so \eqref{eq:TnI=CnI} follows.

If $J=I\setminus \{0\}$, then $J\subseteq[n-1]$.
Let $\sigma'=(r,r+h+1)(\barr{r},\barr{r+h+1})\sigma$. Then $\sigma'\in C_n$ and  $(r,h)$ is also an $a_0$-odd sandwich in $\sigma'$.
It follows from Definition \ref{def:k-oddsandwich} that for all
$$i\in[a_0]\setminus\{|\sigma^{-1}(r)|,|\sigma^{-1}(r+h+1)|\}\quad \text{and}\quad j\in\{|\sigma^{-1}(r)|,|\sigma^{-1}(r+h+1)|\},$$
we have
$$\sigma(i)>\sigma(j)\Leftrightarrow \sigma'(i)>\sigma'(j)\quad \text{and} \quad \sigma(i)+\sigma(j)<0\Leftrightarrow \sigma'(i)+\sigma'(j)<0.$$
Hence $\sigma \in C_n^{J}\setminus T_n^{J}$ if and only if $\sigma' \in C_n^{J}\setminus T_n^{J}$. The $(r,h)$ is chosen to be lexicographically
minimal, so it is unique for $\sigma$.
Thus, $\sigma\mapsto\sigma'$ is a fixed-point free involution on $C_n^{J}\setminus T_n^{J}$.
Since $h$ is odd, it follows that $\oinv(\sigma)=\oinv(\sigma')$ and  $\onsp(\sigma)=\onsp(\sigma')$, so that $L(\sigma')=L(\sigma)$.
But $\ell(\sigma')=\ell(\sigma)\pm 1$, so $(-1)^{\ell(\sigma)}x^{L(\sigma)}+(-1)^{\ell(\sigma')}x^{L(\sigma')}=0$, and the proof follows.
\end{proof}

Let $I=[0,a_0-1]\cup[a_0+1,n-1]$, where $a_0\in[2,n-2]$. We now consider a parabolic factorization  of the form $\sigma=\sigma^{[0,a_0-1]}\sigma_{[0,a_0-1]}$
for each element $\sigma$ in $T_n^{I\setminus\{0\}}$.
It is easy to see that $\sigma_{[0,a_0-1]}\in \mathcal{D}_{a_0}\times \{1_{n-a_0}\}$, where
$$\mathcal{D}_{a_0}\times \{1_{n-a_0}\}=\{\sigma\in \mathcal{D}_{n}\,|\, \sigma(1)\cdots\sigma(a_0)\in \mathcal{D}_{a_0}\ \text{and}\ \sigma(i)=i\ \text{for all}\ i\in[a_0+1,n]\}.$$
Moreover, $\sigma_{[0,a_0-1]}(1)\cdots\sigma_{[0,a_0-1]}(a_0)$ is obtained from $\sigma$ by first taking standardization of $\sigma(1)\cdots\sigma(a_0)$, and then (possibly) changing the signs of $\pm1$ such that the resulting sequence is an even signed permutation. That is, for all $i,j\in[a_0]$ we have
\begin{align}\label{eq:sigma_[0,a_0-1]1}
	|\sigma(i)|>|\sigma(j)|\Leftrightarrow |\sigma_{[0,a_0-1]}(i)|>|\sigma_{[0,a_0-1]}(j)|,
\end{align}
and for all $k\in[a_0]$ with $|\sigma_{[0,a_0-1]}(k)|\neq 1$ we have
\begin{align}\label{eq:sigma_[0,a_0-1]2}
	\sgn(\sigma(k))=\sgn(\sigma_{[0,a_0-1]}(k)).
\end{align}
For example, if $a_0=3$, then $\bar5\bar21\bar4\bar3=125\bar4\bar3\cdot \bar3\bar2145$ and $\bar5\bar2\bar14\bar3=\bar1254\bar3\cdot \bar3\bar2145$.

\begin{lemma}\label{lem:TnIfactorTnI0decom}
	Let $I\subseteq[0,n-1]$ with $I_0=[0,a_0-1]$, where $a_0\geq 3$. If $a_0\equiv n\equiv 1(\bmod\ 2)$, then
\begin{align*}
T_n^I=\{\sigma^{[0,a_0-1]}\,|\, \sigma\in T_n^{I\setminus\{0\}}\},\quad C_{a_0}^{[a_0-1]}\times \{1_{n-a_0}\}=\{\sigma_{[0,a_0-1]}\,|\, \sigma\in T_n^{I\setminus\{0\}}\},
\end{align*}
and $T_n^{I\setminus\{0\}}= T_n^I\cdot (C_{a_0}^{[a_0-1]}\times \{1_{n-a_0}\})$.
\end{lemma}
\begin{proof}
Let $\sigma\in T_n^{I\setminus\{0\}}$. The approach we use to obtain $\sigma^{[0,a_0-1]}$ shows that $\sigma^{[0,a_0-1]}\in \mathcal{D}_n^I$ without $a_0$-odd sandwiches.
Since $n\equiv 1(\bmod\ 2)$, we have $C_n=C_{n,0}$. Let $|\sigma(i_1)|<|\sigma(i_2)|<\cdots<|\sigma(i_{a_{0}})|$ such that
$$\{|\sigma(i_1)|,|\sigma(i_2)|,\ldots,|\sigma(i_{a_{0}})|\}=\{|\sigma(1)|,|\sigma(2)|,\ldots,|\sigma(a_{0})|\}.$$
We claim that there does not exist $k\in[a_0-1]$ such that $|\sigma(i_k)|\equiv|\sigma(i_{k+1})|(\bmod\ 2)$; otherwise,
let $r=|\sigma(i_k)|$ and $h=|\sigma(i_{k+1})|-|\sigma(i_k)|-1$. Then for any $i\in[h]$, we have $r+i\not\in \{|\sigma(j)|\,|\,j\in[a_0]\}$, that is, $|\sigma^{-1}(r+i)|>a_0$.
So $(r,h)$ is an $a_0$-odd sandwich in $\sigma$, contradicting $\sigma\in T_n^{I\setminus\{0\}}$.
Thus, $|\sigma(i_k)|\not\equiv|\sigma(i_{k+1})|(\bmod\ 2)$, and hence $\sigma^{[0,a_0-1]}(k)\not\equiv\sigma^{[0,a_0-1]}(k+1)(\bmod\ 2)$ for all $k\in[a_0-1]$.
But $a_0\equiv 1(\bmod\ 2)$, we see that $\sigma^{[0,a_0-1]}(k)\equiv k(\bmod\ 2)$ for all $k\in[a_0]$.
Since $\sigma^{[0,a_0-1]}(k)= \sigma(k)\equiv k(\bmod\ 2)$ for any $k\in[a_0+1,n]$, we obtain that $\sigma^{[0,a_0-1]}\in C_n$, and hence $\sigma^{[0,a_0-1]}\in T_n^{I}$.
Notice that $C_{n}$ is a subgroup of $\mathcal{D}_n$. So $\sigma_{[0,a_0-1]}\in C_{n}^{[a_0-1]}$, which together with $\sigma_{[0,a_0-1]}\in \mathcal{D}_{a_0}\times \{1_{n-a_0}\}$ yields that $\sigma_{[0,a_0-1]}\in C_{a_0}^{[a_0-1]}\times \{1_{n-a_0}\}$.
Conversely, for any $\tau\in T_n^I$ and $\delta\in C_{a_0}^{[a_0-1]}\times \{1_{n-a_0}\}$, we have $\tau\delta\in  T_n^{I\setminus\{0\}}$ satisfying
$\tau=(\tau\delta)^{[0,a_0-1]}$ and $\delta=(\tau\delta)_{[0,a_0-1]}$, so the proof follows.
\end{proof}

\begin{lemma}\label{lem:TnIfactorL=L+L}
	Let $I\subseteq[0,n-1]$ with $I_0=[0,a_0-1]$, where $a_0\geq 3$. If $a_0\equiv n\equiv 1(\bmod\ 2)$, then for any $\sigma\in T_n^{I\setminus\{0\}}$, we have
	$L(\sigma)=L(\sigma^{[0,a_0-1]})+L(\sigma_{[0,a_0-1]})$.
\end{lemma}
\begin{proof}
	Let $\sigma\in T_n^{I\setminus\{0\}}$ and $\tau=\sigma^{[0,a_0-1]}$. Define
	\begin{align*}
		P&\,=\{(i,j)\in[a_0]\times[a_0+1,n]\,|\,\sigma(i)>\sigma(j),i\not\equiv j(\bmod\ 2)\},\\
		Q&\,=\{(i,j)\in[a_0]\times[a_0+1,n]\,|\,\sigma(i)+\sigma(j)<0,i\not\equiv j(\bmod\ 2)\},\\
		X&\,=\{(i,j)\in[a_0]\times[a_0+1,n]\,|\,\tau(i)>\tau(j),i\not\equiv j(\bmod\ 2)\},\\
		Y&\,=\{(i,j)\in[a_0]\times[a_0+1,n]\,|\,\tau(i)+\tau(j)<0,i\not\equiv j(\bmod\ 2)\}.
	\end{align*}
	Notice that $\tau(i)=\sigma(i)$ for all $i\in[a_0+1,n]$.
	We have
	\begin{align*}
		L(\sigma)=&\,\oinv(\sigma)+\onsp(\sigma)\\
		=&\, |P|+| Q|+|\{(i,j)\in[a_0+1,n]^2\,|\,i<j,\sigma(i)>\sigma(j),i\not\equiv j(\bmod\ 2)\}|\\
		&\,+|\{(i,j)\in[a_0]^2\,|\,i<j,\sigma(i)+\sigma(j)<0,i\not\equiv j(\bmod\ 2)\}|\\
		&\,+|\{(i,j)\in[a_0+1,n]^2\,|\,i<j,\sigma(i)+\sigma(j)<0,i\not\equiv j(\bmod\ 2)\}|,\\
		L(\sigma^{[0,a_0-1]}) =&\,\oinv(\tau)+\onsp(\tau)\\
		=&\,|X|+|Y|+ |\{(i,j)\in[a_0+1,n]^2\,|\,i<j,\sigma(i)>\sigma(j),i\not\equiv j(\bmod\ 2)\}|\\
		&\,+|\{(i,j)\in[a_0+1,n]^2\,|\,i<j,\sigma(i)+\sigma(j)<0,i\not\equiv j(\bmod\ 2)\}|,
	\end{align*}
	and, by \eqref{eq:sigma_[0,a_0-1]1} and \eqref{eq:sigma_[0,a_0-1]2},
	\begin{align*}
		L(\sigma_{[0,a_0-1]})&\ =\oinv(\sigma_{[0,a_0-1]})+\onsp(\sigma_{[0,a_0-1]})\\
		&\ =|\{(i,j)\in[a_0]^2\,|\,i<j,\sigma_{[0,a_0-1]}(i)+\sigma_{[0,a_0-1]}(j)<0,i\not\equiv j(\bmod\ 2)\}|\\
		&\ =|\{(i,j)\in[a_0]^2\,|\,i<j,\sigma(i)+\sigma(j)<0,i\not\equiv j(\bmod\ 2)\}|.
	\end{align*}
	Thus, $L(\sigma)-\,L(\sigma^{[0,a_0-1]})-L(\sigma_{[0,a_0-1]})=|P|+| Q|-| X|-|Y|$.
	Define a map
	$$
	\varphi:P\cup Q \rightarrow X\cup Y,\quad \varphi(i,j)=(|\tau^{-1}(\sigma(i))|,j).
	$$
	By Lemma \ref{lem:TnIfactorTnI0decom}, $\sigma,\tau\in C_n$.
	Since $n\equiv 1(\bmod\ 2)$, we see that $C_n=C_{n,0}$, so that $\tau^{-1}(\sigma(i))\equiv \sigma(i)\equiv i(\bmod\,2)$ for all $i\in[n]$.
	Hence, we conclude from $i\not\equiv j(\bmod\ 2)$ that $|\tau^{-1}(\sigma(i))|\not\equiv j(\bmod\ 2)$ for all $(i,j)\in[n]^2$.
	Note that for all $(i,j)\in P\cup Q$, we have $\tau(j)=\sigma(j)$.
	For any $(i,j)\in P\backslash Q $, we have $\sigma(i)>\sigma(j)>\barr{\sigma(i)}$, so that
	$\varphi(i,j)=(\tau^{-1}(\sigma(i)),j)\in X\backslash Y$ if $\tau^{-1}(\sigma(i))>0$, and
	$\varphi(i,j)=(\tau^{-1}(\overline{\sigma(i)}),j)\in Y\backslash X$ if $\tau^{-1}(\sigma(i))<0$, that is, $\varphi(P\backslash Q)\subseteq  X\backslash Y\cup Y\backslash X$.
	It is clear by symmetry that $\varphi(Q\backslash P)\subseteq  X\backslash Y\cup Y\backslash X$.
	For any $(i,j)\in P\cap Q$, we have  $\sigma(j)<\overline{|\sigma(i)|}$ with $i\not\equiv j(\bmod\ 2)$, so that $\varphi(i,j)\in X\cap Y$.
	Thus, $\varphi$ is well defined.  It is clear that $\varphi$ is a bijection, so $L(\sigma)=L(\sigma^{[0,a_0-1]})+L(\sigma_{[0,a_0-1]})$.
\end{proof}

\begin{prop}\label{prop:noddIcompressed}
Let $I=[0,a_0-1]\cup[a_0+1,n-1]$ with $a_0\in[2,n-2]$. If $a_0\equiv n\equiv 1(\bmod\ 2)$, then
%$\mathcal{D}_n^{I\setminus\{0\}}(x)=\mathcal{D}_n^I(x)  \mathcal{D}_{a_0}^{[a_0-1]}(x)$, so that
\begin{align*}
\mathcal{D}_n^{I}(x)=\left[\begin{matrix}\frac{n-1}{2}\\ \frac{a_0-1}{2},\frac{n-a_0}{2}\end{matrix}\right]_{x^2}\prod_{i=a_0+1}^{n}\left(1+(-1)^{i-1}x^{\left\lfloor\frac{i}{2}\right\rfloor}\right).
\end{align*}
\end{prop}
\begin{proof}
	By Lemma \ref{lem:eqDnI=TnI}, $\mathcal{D}_n^{I\setminus\{0\}}(x)=\sum_{\sigma\in T_n^{I\setminus\{0\}}}(-1)^{\ell(\sigma)}x^{L(\sigma)}$.
	It follows from Lemmas \ref{lem:TnIfactorTnI0decom} and \ref{lem:TnIfactorL=L+L} that
\begin{align*}
\mathcal{D}_n^{I\setminus\{0\}}(x)=\Big(\sum_{\sigma\in T_n^I}(-1)^{\ell(\sigma)}x^{L(\sigma)}\Big)\Big(\sum_{\sigma\in C_{a_0}^{[a_0-1]}\times\{1_{n-a_0}\}}(-1)^{\ell(\sigma)}x^{L(\sigma)}\Big).
\end{align*}
Hence, by Lemmas \ref{lem:DnI=CnIeq} and \ref{lem:eqDnI=TnI},
$\mathcal{D}_n^{I\setminus\{0\}}(x)=\mathcal{D}_n^I(x)  \mathcal{D}_{a_0}^{[a_0-1]}(x)$, so the desired identity follows from Proposition \ref{prop:DnI(x)-Iin[n-1]}.
\end{proof}

\subsection{When $|I_0|$ is an even positive integer}\label{I0=evenpositiveint}
In this subsection, we establish a closed product formula of $\mathcal{D}_n^I(x)$ in the case
where $|I_0|\geq2$ is an even positive integer.

\begin{lemma}\label{lem:DnJ(x)-2components}
	Let $I=[0,a_0-1]\cup[a_0+1,n-1]$, where $a_0\in [2,n-2]$. If $a_0\equiv n\equiv 0(\bmod\ 2)$, then
$$
	\mathcal{D}_n^I(x)=\ x^{n-a_0} \mathcal{D}_{n-1}^{[0,a_0-2]\cup[a_0,n-2]}(x)+(1+x^{a_0}-2x^{a_0+\frac{n}{2}}) \mathcal{D}_{n-1}^{[0,a_0]\cup[a_0+2,n-2]}(x).
$$
\end{lemma}
\begin{proof}
Since $a_0\in [2,n-2]$, we see that $|\sigma(1)|<\sigma(2)$ and hence $\sigma(1)\neq \pm n$. By analyzing the position of $\pm n$ in the one-line notation of $\sigma$, we see that
$\mathcal{D}_n^I(x)$ is equal to
\begin{align}\label{eq:D_n^J-partition-by-n}
	\sum_{\substack{\sigma\in \mathcal{D}_n^I,\, \sigma(a_0)=n}}(-1)^{\ell(\sigma)}x^{L(\sigma)}+\sum_{\substack{\sigma\in \mathcal{D}_n^I,\,
\sigma(a_0+1)=\overline n}}(-1)^{\ell(\sigma)}x^{L(\sigma)}
	+\sum_{\substack{\sigma\in \mathcal{D}_n^I,\, \sigma(n)=n}}(-1)^{\ell(\sigma)}x^{L(\sigma)}.
\end{align}
We next compute the three summands in \eqref{eq:D_n^J-partition-by-n} individually.

(i) Given $\sigma\in \mathcal{D}_n^I$, define $\sigma'=(1,\overline 1)(n,\overline{n})\sigma$. If $a_0=2$, then $I=[0,1]\cup[3,n-1]$, so
$\sigma\in \mathcal{D}_n^I$ with $\sigma(a_0)=n$ is equivalent to
$\sigma'\in \mathcal{D}_n^{[2,n-1]}$ with $\sigma'(a_0)=\overline{n}$. Since $a_0\equiv n(\bmod\ 2)$, it follows from Lemmas \ref{lem:ellLsigma'a=n} and \ref{lem:ai=oddn}\eqref{item:sigmaa=-na} that
\begin{align*}
	\sum_{\substack{\sigma\in \mathcal{D}_n^I\\ \sigma(a_0)=n}}(-1)^{\ell(\sigma)}x^{L(\sigma)}
	= \frac{1}{x^{2}} \sum_{\substack{\sigma'\in \mathcal{D}_n^{[2,n-1]}\\ \sigma'(2)=\overline{n}}}(-1)^{\ell(\sigma')}x^{L(\sigma')}
    =  \frac{1}{x^{2}} \sum_{\substack{\sigma\in \mathcal{D}_n^{[2,n-1]}\\ \sigma(2)=\overline{n},\, \sigma(n)<\sigma(1)}}(-1)^{\ell(\sigma)}x^{L(\sigma)}.
\end{align*}
For any $\sigma\in \mathcal{D}_{n-1}^{[n-2]}$, let $\sigma''=(1,\bar1)\cdot \left(\sigma(n-1)\overline{n}\sigma(1)\cdots\sigma(n-2)\right)$, that is, $\sigma''$ is obtained from $\sigma(n-1)\overline{n}\sigma(1)\cdots\sigma(n-2)$ by interchanging the signs of $\pm 1$.
Then $\sigma\in \mathcal{D}_{n-1}^{[n-2]}$ if and only if $\sigma''\in \mathcal{D}_{n}^{[2,n-1]}$ with $\sigma''(2)=\overline{n}$ and $\sigma''(n)<\sigma''(1)$.
Since $\ell(\sigma'')=\ell(\sigma)+2n-2$ and
$L(\sigma'')=L(\sigma)+n$, it follows from  Lemma \ref{prop:D_n^Icup0=Icup1} that
\begin{align*}
\sum_{\substack{\sigma\in \mathcal{D}_n^I,\, \sigma(a_0)=n}}(-1)^{\ell(\sigma)}x^{L(\sigma)}
	= x^{n-2} \mathcal{D}_{n-1}^{[n-2]}(x)
=x^{n-a_0}\mathcal{D}_{n-1}^{[0,a_0-2]\cup[a_0,n-2]}(x).
\end{align*}
If $a_0\geq 4$, then $\sigma\in \mathcal{D}_n^I$ with $\sigma(a_0)=n$ if and only if  $\sigma'\in \mathcal{D}_n^{[0,a_0-2]\cup[a_0,n-1]}$ with $\sigma'(a_0)=\overline{n}$,
where $\sigma'=(1,\overline 1)(n,\overline{n})\sigma$.
Hence,
\begin{align*}
	\sum_{\substack{\sigma\in \mathcal{D}_n^I\\ \sigma(a_0)=n}}(-1)^{\ell(\sigma)}x^{L(\sigma)}
	=&\, \frac{1}{x^{a_0}} \sum_{\substack{\sigma'\in \mathcal{D}_n^{[0,a_0-2]\cup[a_0,n-1]}\\ \sigma'(a_0)=\overline{n}}}(-1)^{\ell(\sigma')}x^{L(\sigma')}
&&{\rm by\ Lemma\ \ref{lem:ellLsigma'a=n}}\\
	=&\, \frac{1}{x^{a_0}(1+x^{a_0-2})} \sum_{\substack{\sigma\in \mathcal{D}_n^{[0,a_0-3]\cup[a_0,n-1]}\\ \sigma(a_0)=\overline{n}}}(-1)^{\ell(\sigma)}x^{L(\sigma)}
&&{\rm by\ Lemma\ \ref{lem:a0=even}\eqref{itemeq:DnxIa0=nb}}\\
=&\, \frac{1}{x^{a_0}(1+x^{a_0-2})}\sum_{\substack{\sigma\in \mathcal{D}_n^{[0,a_0-3]\cup[a_0,n-1]}\\ \sigma(a_0)=\overline{n},\, \sigma(n)<\sigma(a_0-1)}}(-1)^{\ell(\sigma)}x^{L(\sigma)}
&&{\rm by\ Lemma\ \ref{lem:ai=oddn}\eqref{item:sigmaa=-na}}\\
=&\, \frac{1}{x^{a_0}(1+x^{a_0-2})}\sum_{\substack{\sigma\in \mathcal{D}_{n-1}^{[0,a_0-3]\cup[a_0-1,n-2]}}}(-1)^{\ell(\sigma''')}x^{L(\sigma''')},
\end{align*}
where $$\sigma'''=(1,\bar1)\cdot \left(\sigma(1)\cdots\sigma(a_0-2)\sigma(n-1)\overline{n}\sigma(a_0-1)\cdots\sigma(n-2)\right)$$ for any $\sigma\in \mathcal{D}_{n-1}^{[0,a_0-3]\cup[a_0-1,n-2]}$. Now
$\ell(\sigma''')=\ell(\sigma)+2n-2$ and $L(\sigma''')= L(\sigma)+n$, from which we see that
\begin{align*}
\sum_{\substack{\sigma\in \mathcal{D}_{n-1}^{[0,a_0-3]\cup[a_0-1,n-2]}}}(-1)^{\ell(\sigma''')}x^{L(\sigma''')}
=&x^n \mathcal{D}_{n-1}^{[0,a_0-3]\cup[a_0-1,n-2]}(x)\\
=&x^n(1+x^{a_0-2})\mathcal{D}_{n-1}^{[0,a_0-2]\cup[a_0,n-2]}(x),
\end{align*}
where the second equality follows from Proposition \ref{prop:a0=even}. Thus,
\begin{align}\label{eq:D0a0-1J1}
	\sum_{\substack{\sigma\in \mathcal{D}_n^I,\, \sigma(a_0)=n}}(-1)^{\ell(\sigma)}x^{L(\sigma)}
	= x^{n-a_0}\mathcal{D}_{n-1}^{[0,a_0-2]\cup[a_0,n-2]}(x).
\end{align}
So equation \eqref{eq:D0a0-1J1} is true for all  $a_0\in [2,n-2]$ with $a_0\equiv n\equiv 0(\bmod\ 2)$.

(ii) Recall that for each $\sigma\in \mathcal{D}_n$, we have defined $\sigma'=(1,\bar 1)(n,\overline{n})\sigma$. Clearly, $\sigma\in \mathcal{D}_n^I$ with $\sigma(a_0+1)=\overline{n}$
if and only if
$\sigma'\in \mathcal{D}_n^{[0,a_0]\cup[a_0+2,n-1]}$ with $\sigma'(a_0+1)=n$.
Since $a_0\equiv n\equiv 0(\bmod\ 2)$, there follows
\begin{align}\label{eq:D0a0-1J2}
	\sum_{\substack{\sigma\in \mathcal{D}_n^I\\ \sigma(a_0+1)=\overline{n}}}(-1)^{\ell(\sigma)}x^{L(\sigma)}
&= x^{a_0} \sum_{\substack{\sigma\in \mathcal{D}_n^{[0,a_0]\cup[a_0+2,n-1]}\\ \sigma(a_0+1)=n}}(-1)^{\ell(\sigma)}x^{L(\sigma)}&&{\rm by\ Lemma\ \ref{lem:ellLsigma'a=n}}\notag\\
&=x^{a_0} \sum_{\substack{\sigma\in \mathcal{D}_n^{[0,a_0]\cup[a_0+2,n-2]}\\ \sigma(a_0+1)=n}}(-1)^{\ell(\sigma)}x^{L(\sigma)}&&
{\rm by\ Lemma\ \ref{lem:tocompressset}\eqref{item:comprs2}}\notag\\
&=-2x^{a_0+\frac{n}{2}} \mathcal{D}_{n-1}^{[0,a_0]\cup[a_0+2,n-2]}(x).&& {\rm by\ Lemma\ \ref{lem:a0=oddn}}
\end{align}

(iii) For the third term of \eqref{eq:D_n^J-partition-by-n}, it follows from Lemmas \ref{prop:shift[0,n-1]}\eqref{eq:prop3.7dnx0} and \ref{lem:a0=even}\eqref{itemeq:DnxIa0a} that
\begin{align}\label{eq:D0a0-1J3}
	\sum_{\substack{\sigma\in \mathcal{D}_n^I,\, \sigma(n)=n}}(-1)^{\ell(\sigma)}x^{L(\sigma)}
= \mathcal{D}_{n-1}^{[0,a_0-1]\cup[a_0+1,n-2]}(x)=(1+x^{a_0})\mathcal{D}_{n-1}^{[0,a_0]\cup[a_0+2,n-2]}(x).
\end{align}
Substituting \eqref{eq:D0a0-1J1}, \eqref{eq:D0a0-1J2} and \eqref{eq:D0a0-1J3} into \eqref{eq:D_n^J-partition-by-n} yields the desired identity.
\end{proof}

\begin{prop}\label{prop:nevenIcompressed}
	Let $I=[0,a_0-1]\cup[a_0+1,n-1]$, where $a_0\in [2,n-2]$. If $a_0\equiv n\equiv 0(\bmod\ 2)$, then
\begin{align*}
	\mathcal{D}_n^I(x)=C_I
\frac{1+x^{a_0}+2x^{\frac{n}{2}}}{1+x^{\frac{n}{2}}}
\prod_{i=a_0+2}^{n} \left(1+(-1)^{i-1}x^{\left\lfloor\frac{i}{2}\right\rfloor}\right).
\end{align*}
\end{prop}
\begin{proof}
Since $a_0\equiv n\equiv0(\bmod\,2)$, it follows from  Lemmas \ref{prop:D_n^Icup0=Icup1} and \ref{prop:Dn[n-1](x)} that if $a_0=2$, then
$$
\mathcal{D}_{n-1}^{[0,a_0-2]\cup[a_0,n-2]}(x)=\mathcal{D}_{n-1}^{[n-2]}(x)=\prod_{i=2}^{n-1}\left(1+(-1)^{i-1}x^{\left\lfloor\frac{i}{2}\right\rfloor}\right),
$$
and from Proposition \ref{prop:noddIcompressed} we conclude that if $a_0\geq4$, then
\begin{align}\label{eq:Dn-10a0-2n-2x}
\mathcal{D}_{n-1}^{[0,a_0-2]\cup[a_0,n-2]}(x)
=&\left[\begin{matrix}\frac{n-2}{2}\\ \frac{a_0-2}{2},\frac{n-a_0}{2}\end{matrix}\right]_{x^2}
\prod_{i=a_0}^{n-1} \left(1+(-1)^{i-1}x^{\left\lfloor\frac{i}{2}\right\rfloor}\right).
\end{align}
So \eqref{eq:Dn-10a0-2n-2x} is true for any $a_0\geq2$.
Substituting $a_0+2$ for $a_0$ in \eqref{eq:Dn-10a0-2n-2x} gives
\begin{align*}
 \mathcal{D}_{n-1}^{[0,a_0]\cup[a_0+2,n-2]}(x)=&\left[\begin{matrix}\frac{n-2}{2}\\ \frac{a_0}{2},\frac{n-a_0-2}{2}\end{matrix}\right]_{x^2}
\prod_{i=a_0+2}^{n-1} \left(1+(-1)^{i-1}x^{\left\lfloor\frac{i}{2}\right\rfloor}\right).
\end{align*}
\delete{Notice that
\begin{align*}
\left[\begin{matrix}\frac{n-2}{2}\\ \frac{a_0-2}{2},\frac{n-a_0}{2}\end{matrix}\right]_{x^2}=C_I\frac{1-x^{a_0}}{1-x^{n}},\quad
\left[\begin{matrix}\frac{n-2}{2}\\ \frac{a_0}{2},\frac{n-a_0-2}{2}\end{matrix}\right]_{x^2}=C_I\frac{1-x^{n-a_0}}{1-x^{n}}.
\end{align*}}
The proposition now follows from Lemma \ref{lem:DnJ(x)-2components}.
\end{proof}

\subsection{Closed product formulas for $\mathcal{D}_n^I(x)$}\label{mainthm}

Now we are in a position to give our main result. Since $\mathcal{D}_n^{[0,n-1]}(x)=1$ and
by Lemma \ref{prop:Dn[n-1](x)} and \cite[Theorem 4.1]{BC17D},
$\mathcal{D}_n^{\emptyset}(x)=\mathcal{D}_n^{[n-1]}(x)^2$,
we assume that $I$ is a proper subset of $[0,n-1]$ in this subsection.

%=================第一种形式=====================================================
\begin{theorem}\label{theo:finaltheorem}
Let $I$ be a proper subset of $[0,n-1]$, where $n\geq2$.
Denote that $\widetilde I=I$ if $I_0=\emptyset$, and $\widetilde I=(I\backslash \{0\})\cup\{1\}$ otherwise. Then
\begin{align}\label{eq:D_nIx2inII_0=1}
	&D_{n}^I(x)=f(x)C_{\widetilde{I}}\prod\limits_{i=2\left\lfloor\frac{|I_0|+2}{2}\right\rfloor}^{n} \left(1+(-1)^{i-1}x^{\left\lfloor\frac{i}{2}\right\rfloor}\right)\prod\limits_{j=2m_{\widetilde{I}}+2}^n\left(1+(-1)^{j-1}x^{\left\lfloor\frac{j}{2}\right\rfloor}\right),
\end{align}
where
\begin{align*}
f(x)&=\begin{cases}
\frac{1+x^{|I_0|}+2x^{m_I}}{1+x^{m_I}}, &{\rm if}\ |I_0|\geq2,\ |I_0|\equiv 0(\bmod\ 2),\ {\rm and}\ n=2m_I,\\ \vspace{1mm}
1+x^{|I_0|}, &{\rm if}\ |I_0|\geq2,\ |I_0|\equiv 0(\bmod\ 2),\ {\rm and}\ n>2m_I,\\ \vspace{1mm}
1, &{\rm otherwise}.
\end{cases}
\end{align*}
\end{theorem}
\begin{proof}
If $|I_0|=0$, that is, $I_0=\emptyset$, then $I\subseteq[n-1]$ and the desired identity follows from Proposition \ref{prop:DnI(x)-Iin[n-1]}.
If $|I_0|=1$, that is, $I_0=\{0\}$, then, by Lemma \ref{prop:D_n^Icup0=Icup1}, $\mathcal{D}_n^I(x)=\mathcal{D}_n^{\widetilde I}(x)$ and the proof also follows from Proposition \ref{prop:DnI(x)-Iin[n-1]}.

Now assume that $|I_0|\geq2$. Then $\widetilde{I}=I\backslash \{0\}$.
We write $I'=C(I)=[0,a_0-1]\cup[a_0+1,a_1-1]\cup\cdots\cup[a_{s-1}+1,a_s-1]$. Then $|I'_0|=|I_0|=a_0$
 and $a_0\equiv a_s(\bmod\ 2)$.
Let $J=[0,a_0-1]\cup[a_0+1,a_s-1]$ and $K=[1,a_1-a_0-1]\cup\cdots\cup[a_{s-1}-a_0+1,a_s-a_0-1]$.

If $|I_0|\equiv 1(\bmod\ 2)$, then, by \eqref{eq:m_I} and \eqref{eq:compressm}, $m_{I'}=m_I=\frac{a_s+1}{2}$, $m_{\widetilde I}=\frac{a_s-1}{2}$.
It follows from Propositions \ref{prop:a0=odd} and \ref{lem:DnI=DnJxSn-aKeq} that
\begin{align*}
\mathcal{D}_n^I(x)=&\,\mathcal{D}_{a_s}^{I'}(x)\prod_{j=a_s+1}^n\left(1+(-1)^{j-1}x^{\lfloor\frac{j}{2}\rfloor}\right)^2\\
=&\,\left[\begin{matrix}\frac{a_s-a_0}{2}\\ \frac{a_1-a_0}{2},\ldots,\frac{a_s-a_{s-1}}{2}\end{matrix}\right]_{x^2}\mathcal{D}_{a_s}^{J}(x)\prod_{j=a_s+1}^n\left(1+(-1)^{j-1}x^{\lfloor\frac{j}{2}\rfloor}\right)^2.
\end{align*}
Applying Proposition \ref{prop:noddIcompressed} gives
%\begin{align*}
%\mathcal{D}_{a_s}^{J}(x)=\left[\begin{matrix}\frac{a_s-1}{2}\\ %\frac{a_0-1}{2},\frac{a_s-a_0}{2}\end{matrix}\right]_{x^2}\prod_{j=a_0+1}^{a_s}\left(1+(-1)^{j-1}x^{\lfloor\frac{j}{2}\rfloor}\right),
%\end{align*}
\eqref{eq:D_nIx2inII_0=1} in the case $|I_0|\equiv 1(\bmod\ 2)$.

The remaining case is that $|I_0|\geq 2$ with $|I_0|\equiv 0(\bmod\ 2)$. It is clear that $m_{I'}=m_{\widetilde{I}}=m_I=\frac{a_s}{2}$. If $n=2m_I$, then $n=a_s$.
It follows from Lemma \ref{lem:tocompressset}\eqref{item:comprs1} that $\mathcal{D}_n^I(x)=\mathcal{D}_{a_s}^{I'}(x)$.
 Then, by Propositions \ref{lem:DnI=DnJxSn-aKeq} and \ref{prop:nevenIcompressed},
\begin{align*}
	\mathcal{D}_n^I(x)
=\left[\begin{matrix}\frac{n-a_0}{2}\\ \frac{a_1-a_0}{2},\ldots,\frac{n-a_{s-1}}{2}\end{matrix}\right]_{x^2}\mathcal{D}_{a_s}^{J}(x)
=C_{\widetilde{I}}\frac{1+x^{a_0}+2x^{\frac{a_s}{2}}}{1+x^{\frac{a_s}{2}}}
\prod_{i=a_0+2}^{n} \left(1+(-1)^{i-1}x^{\left\lfloor\frac{i}{2}\right\rfloor}\right).
\end{align*}
If $n>2m_I$, then $a_s<n$. From Proposition \ref{prop:a0=even}, we deduce that
\begin{align}\label{maineq:as<nmod02}
\mathcal{D}_n^I(x)=\left(1+x^{a_0}\right)\mathcal{D}_{2m_{I''}-1}^{I''}(x)\prod_{j=2m_{I''}}^n\left(1+(-1)^{j-1}x^{\left\lfloor\frac{j}{2}\right\rfloor}\right)^2,
\end{align}
where $I''=[0,a_0]\cup[a_0+2,a_1]\cup\cdots\cup[a_{s-1}+2,a_s]$.
Since $|I''_0|\equiv 1(\bmod\ 2)$, replacing $I$ by $I''$ in \eqref{eq:D_nIx2inII_0=1} gives
\begin{align*}
\mathcal{D}_{2m_{I''}-1}^{I''}(x)=C_{I''\backslash\{0\}}\prod_{j=a_0+2}^{a_s+1}\left(1+(-1)^{j-1}x^{\left\lfloor\frac{j}{2}\right\rfloor}\right).
\end{align*}
Substituting into \eqref{maineq:as<nmod02} gives \eqref{eq:D_nIx2inII_0=1} for the case where  $n>2m_I$,
completing the proof.
\end{proof}

%In \cite{BC17D}, Brenti and Carnevale conjectured several closed product formulae of $\mathcal{D}_n^I(x)$ when $|I\backslash\{0,1\}|\leq 1$.
Now we verify the conjectures of Brenti and Carnevale as consequences of Theorem \ref{theo:finaltheorem}.
\begin{coro}
Conjecture \ref{conj:DnI0i} is true.
\end{coro}
\begin{proof}
Let $\widetilde{I}=\{1,i\}$. Since $i\geq3$, we have $m_{\widetilde{I}}=2$ and $C_{\widetilde{I}}=1+x^2.$
By Theorem \ref{theo:finaltheorem},
\begin{align*}
\mathcal{D}_n^{\{0,i\}}(x)=C_{\widetilde{I}}
\prod_{i=2}^{n} \left(1+(-1)^{i-1}x^{\left\lfloor\frac{i}{2}\right\rfloor}\right)\prod_{j=6}^n \left(1+(-1)^{j-1}x^{\left\lfloor\frac{j}{2}\right\rfloor}\right)
= \prod_{j=4}^n\left(1+(-1)^{j-1}x^{\left\lfloor\frac{j}{2}\right\rfloor}\right)^2
\end{align*}
and
\begin{align*}
\mathcal{D}_n^{\{0,1,i\}}(x)=&\,(1+x^2)C_{\widetilde{I}}
\prod_{j=4}^n\left(1+(-1)^{j-1}x^{\left\lfloor\frac{j}{2}\right\rfloor}\right)\prod_{j=6}^n\left(1+(-1)^{j-1}x^{\left\lfloor\frac{j}{2}\right\rfloor}\right)\\
=&\,(1-x^4)\prod_{j=5}^n\left(1+(-1)^{j-1}x^{\left\lfloor\frac{j}{2}\right\rfloor}\right)^2,
\end{align*}
as required.
\end{proof}

It was predicted in Conjecture \ref{conj:DnI(x)} that $\prod_{j=2m_I+2}^n \left(1+(-1)^{j-1}x^{\left\lfloor\frac{j}{2}\right\rfloor}\right)^2$ is always a divisor of $\mathcal{D}_n^I(x)$.
This is a straightforward consequence of Theorem \ref{theo:finaltheorem}.

\begin{coro}\label{coro:formDnI(x)}
Let $I\subseteq [0,n-1]$ with $n\geq 3$. Then there exists a polynomial $M_I(x)\in \mathbb{Z}[x]$ such that
\begin{align*}
\mathcal{D}_n^I(x)=M_I(x)\prod_{j=2m_I+2}^n \left(1+(-1)^{j-1}x^{\left\lfloor\frac{j}{2}\right\rfloor}\right)^2.
\end{align*}
\end{coro}
\begin{proof}
It suffices to assume that $n\geq 2m_I+2$, since $\prod_{j=2m_I+2}^n \left(1+(-1)^{j-1}x^{\left\lfloor\frac{j}{2}\right\rfloor}\right)=1$ if $n\leq 2m_I+1$. The proof follows directly from Theorem \ref{theo:finaltheorem}.
\end{proof}

\begin{remark}
Let $n\geq 3$ and $I\subseteq [0,n-1]$ with connected components $I_0, I_1,\ldots, I_s$. If $I_0\neq\{0\}$, then, by Theorem \ref{theo:finaltheorem},
the polynomial $M_I(x)$ only depends on $(|I_0|,|I_1|,\ldots,|I_s|)$ and is a symmetric function of $|I_1|,\ldots,|I_s|$.
However, these can fail if $I_0=\{0\}$, as illustrated by the following examples.

Let $n=4$, $I=\{0,2\}$, and $I'=\{0,3\}$. Then $(|I_0|,|I_1|)=(|I'_0|,|I'_1|)=(1,1)$. By Theorem \ref{theo:finaltheorem}, we have
$M_I(x)=(1-x^2)^3$ and $M_{I'}(x)=(1-x^2)^2(1+x^2)$.
So in general $M_I(x)$ does not only depend on $(|I_0|,|I_1|,\ldots,|I_s|)$.

Let $n=6$, $K=\{0,2,3,5\}$, and $K'=\{0,2,4,5\}$. Then
$(|K_0|,|K_1|,|K_2|)=(1,2,1)$ and $(|K'_0|,|K'_1|,|K'_2|)=(1,1,2)$.
Then
$M_K(x)=(1-x^3)(1-x^4)(1-x^6)$
and
$M_{K'}(x)=(1-x^3)^2(1-x^4)^2$,
which indicates that $M_I(x)$ is in general not a symmetric function of $|I_1|,\ldots,|I_s|$.
However, by Theorem \ref{theo:finaltheorem}, both $M_I(x)$ and $\mathcal{D}_n^I(x)$ are symmetric with respect to $|I_2|,\ldots,|I_s|$.
\end{remark}

\section{Products of cyclotomic polynomials}\label{sec:cyclotomicpoly}

 Recall that for any positive integer $n$, the \emph{$n$-th cyclotomic polynomial} \cite{Hun74} is the monic polynomial
$(x-\omega_1)(x-\omega_2)\cdots(x-\omega_r)$,
where $\omega_1,\omega_2,\ldots,\omega_r$ are all the distinct primitive $n$-th roots of unity.
In this section, we give a combinatorial description of the sign-twisted generating functions in type $D$ that are not expressible as products of cyclotomic polynomials, and thereby verify
Conjecture \ref{conj:Ste197.4a}.

\begin{lemma}\label{lem:cycltomicpolyqua}
Let $f(x)=x^n+2x^m+1$, where $m$ and $n$ are positive integers. Then $f(x)$ is a product of cyclotomic polynomials if and only if $n=2m$.
\end{lemma}
\begin{proof}
If $n=2m$, then $f(x)=(x^m+1)^2$, which is clearly a product of cyclotomic polynomials. Conversely,
suppose $f(x)$ is a product of cyclotomic polynomials. Let $d$ be the positive greatest common divisor of $m$ and $n$, and let $m_1=\frac{m}{d}$ and $n_1=\frac{n}{d}$.
Let $t=x^d$, and write $g(t)=f(x)=t^{n_1}+2t^{m_1}+1$.
Take an arbitrary complex root, say $\alpha$, of $g(t)$ and assume that $\omega$ is a $d$-th root of $\alpha$.
Then $f(\omega)=g(\alpha)=0$, so that $-2\omega^m=\omega^n+1$. But $f(x)$ is a product of cyclotomic polynomials, we have $|\omega|=1$, and hence $|\omega^n+1|=|-2\omega^m|=2$. Thus, $\omega^n=1$, so that $\omega^m=-1$.
Since $d$ is the positive greatest common divisor of $m$ and $n$, there exist integers $u$ and $v$ such that $d=um+vn$, so that $\omega^d=\omega^{um+vn}=(-1)^u$,
which together with $(\omega^d)^{m_1}=\omega^m=-1$ yields that $\omega^d=-1$. Thus, $\alpha=-1$, in other words, $g(t)$ has a unique complex root $-1$.
Since $f(x)$ is monic, we have $m<n$, it follows that $g(t)$ is a polynomial of degree $n_1$, so that $g(t)=(t+1)^{n_1}$.
But there are only three non-zero terms in the expansion of $g(t)$, we must have $n_1=2$. Hence $g(t)=t^{2}+2t+1$, from which we see that $m_1=1$. Thus, $n=2m$, completing the proof.
\end{proof}

\begin{theorem}\label{thm:cnanotbeprodcyclopoly}
Let $I$ be a proper subset of $[0,n-1]$, where $n\geq2$. Then $\mathcal{D}_n^I(x)$ can not be expressed as a product of cyclotomic polynomials if and only if $n\equiv0(\bmod\ 2)$ and $\{0\}\cup\{i\in[n-1]\,|\, i\equiv1(\bmod\ 2)\}\subseteq I$.
In this case, we have
$$\mathcal{D}_n^I(x)=C_I\frac{1+x^{|I_0|}+2x^{\frac{n}{2}}}{1+x^{\frac{n}{2}}}
\prod\limits_{i=|I_0|+2}^{n} \left(1+(-1)^{i-1}x^{\left\lfloor\frac{i}{2}\right\rfloor}\right).$$
\end{theorem}
\begin{proof}
If $\mathcal{D}_n^I(x)$ can not be expressed as a product of cyclotomic polynomials, then, by Theorem \ref{theo:finaltheorem}, $|I_0|\geq2$, $|I_0|\equiv 0(\bmod\ 2)$, and $n=2m_I$.
Let $I_0=[0,a_0-1]$ and $I_k=[a_k+1,b_k-1]$, where $k\in [s]$. Then $a_0\equiv n\equiv0(\bmod\, 2)$, so we have
$$
m_I=\sum_{k=0}^s\left\lfloor\frac{|I_k|+1}{2}\right\rfloor=\frac{a_0}{2}+\sum_{k=1}^s\left\lfloor\frac{b_k-a_k}{2}\right\rfloor\leq \frac{b_s}{2}\leq \frac{n}{2}.
$$
Since $n=2m_I$, we see that $a_0=a_1$, $b_s=n$, and $a_{k+1}=b_k\equiv 0(\bmod\ 2)$ for all $k\in [s-1]$.
Thus, $\{0\}\cup\{i\in[n-1]\,|\, i\equiv1(\bmod\ 2)\}\subseteq I$.

Conversely, assume that $n\equiv0(\bmod\ 2)$ and $\{0\}\cup\{i\in[n-1]\,|\, i\equiv1(\bmod\ 2)\}\subseteq I$.
Then there exist $a_0,a_1,\ldots,a_{s-1}\in [n-1]$ such that $a_0\equiv a_1\equiv \cdots\equiv a_{s-1}\equiv n(\bmod\ 2)$ and $I=[0,a_0-1]\cup[a_0+1,a_1-1]\cup\cdots\cup[a_{s-1}+1,n-1]$.
By \eqref{eq:m_I}, $m_I=\frac{n}{2}$.
It follows from Theorem \ref{theo:finaltheorem} that
\begin{align*}
\mathcal{D}_n^I(x)&=C_{I}\frac{1+x^{|I_0|}+2x^{\frac{n}{2}}}{1+x^{\frac{n}{2}}}
\prod\limits_{i=|I_0|+2}^{n} \left(1+(-1)^{i-1}x^{\left\lfloor\frac{i}{2}\right\rfloor}\right).
\end{align*}
Since $I\subsetneq [0,n-1]$, we have $s\geq1$, so that $m_I>\frac{|I_0|}{2}$, that is, $|I_0|<2m_I$. By Lemma \ref{lem:cycltomicpolyqua},
$1+x^{|I_0|}+2x^{m_I}$ is not a product of cyclotomic polynomials. Thus, $\mathcal{D}_n^I(x)$ can not be expressed as a product of cyclotomic polynomials.
\end{proof}

\begin{coro}\label{coro:Stembconj}
Conjecture \ref{conj:Ste197.4a} is true.
\end{coro}
\begin{proof}
If $\Phi\cong D_{2m}$ for some positive integer $m$, then $\Phi_K\cong A_1^{\oplus m+1}$ is equivalent to $K=\{0\}\cup\{i\in[2m-1]\,|\, i\equiv1(\bmod\ 2)\}$.
Hence Conjecture \ref{conj:Ste197.4a} follows from Theorem \ref{thm:cnanotbeprodcyclopoly}.
\end{proof}

\noindent
{\bf Acknowledgements.}
%The authors would like to thank Li Guo for helpful conversations.
This work was partially supported by the National Natural Science Foundation of China (Grant Nos. 12471020 and 12071377),
the Natural Science Foundation of Chongqing (Grant No. CSTB2023NSCQ-MSX0706), and the Fundamental Research Funds for the Central Universities (Grant No. SWU-XDJH202305).


\begin{thebibliography}{abcdsfgh}

\bibitem{AGR05} R. Adin, I. Gessel, Y. Roichman, Signed Mahonians, J. Combin. Theory Ser. A 109(2005), 25--43.

\bibitem{Bia06} R. Biagioli, Signed Mahonian polynomials for classical Weyl groups, European J. Combin. 27(2006), 207--217.

\bibitem{BB05} A. Bj\"orner, F. Brenti, Combinatorics of Coxeter Groups, in: Graduate Texts in Mathematics, vol. 231, Springer-Verlag, New York, 2005.

\bibitem{BC17D} F. Brenti, A. Carnevale, Odd length for even hyperoctahedral groups and signed generating functions, Discrete Math. 340(2017), 2822--2833.

\bibitem{BC21} F. Brenti, A. Carnevale, Odd length: Odd diagrams and descent classes, Discrete Math. 344(2021), 112308.

\bibitem{BC19} F. Brenti, A. Carnevale, Odd length in Weyl groups, Algebr. Comb. 2(2019), 1125--1147.

\bibitem{BC17T} F. Brenti, A. Carnevale, Proof of a conjecture of Klopsch-Voll on Weyl groups of type A, Trans. Amer. Math. Soc. 369(2017), 7531--7547.

\bibitem{BCT23}  F. Brenti, A. Carnevale, B.E. Tenner, Odd diagrams, Bruhat order, and pattern avoidance, Combin. Theory 2(2022), \#13.

\bibitem{BC20} F. Brenti, P. Sentinelli,  Odd and even major indices and one-dimensional characters for classical Weyl groups, Ann. Comb. 24(2020), 809--835.

\bibitem{Cas12} F. Caselli, Signed mahonians on some trees and parabolic quotients, J. Combin. Theory Ser. A 119(2012), 1447--1460.

\bibitem{DF92} J. D\'esarm\'enien, D. Foata, The signed Eulerian numbers, Discrete Math. 99(1992), 49--58.

\bibitem{FG23} N.J.Y. Fan, P.L. Guo, Poincar\'e polynomials of odd diagram classes, SIAM J. Discrete Math. 36(2022), 2225--2237.

\bibitem{Hum90} J.E. Humphreys, Reflection Groups and Coxeter Groups, Cambridge University Press, Cambridge, 1990.

\bibitem{Hun74} T.W. Hungerford, Algebra, in: Graduate Texts in Mathematics, vol. 73, Springer-Verlag, New York, 1974.

\bibitem{Igu00} J.-I. Igusa, An Introduction to the Theory of Local Zeta Functions, AMS/IP Studies
in Advanced Mathematics, vol. 14, American Mathematical Society, Providence, RI, 2000.

\bibitem{KV09} B. Klopsch, C. Voll, Igusa-type functions associated to finite formed spaces and their functional equations, Trans. Amer. Math. Soc. 361(2009), 4405--4436.

\bibitem{Lan18} A. Landesman, Proof of Stasinski and Voll's hyperoctahedral group conjecture, Australas. J. Combin. 71(2018), 196--240.

\bibitem{Mon15} P. Mongelli, Signed excedance enumeration in classical and affine Weyl groups, J. Combin. Theory Ser. A 130(2015), 129--149.

\bibitem{Rei95} V. Reiner, Descents and one-dimensional characters for classical Weyl groups, Discrete Math. 140(1995), 129--140.

\bibitem{Siv11} S. Sivasubramanian, Signed excedance enumeration via determinants, Adv. in Appl. Math. 47(2011), 783--794.

\bibitem{SV13} A. Stasinski, C. Voll, A new statistic on the hyperoctahedral groups, Electron. J. Combin. 20(2013), \#P50.

\bibitem{SV14} A. Stasinski, C. Voll, Representation zeta functions of nilpotent groups and generating functions for Weyl groups of type B, Amer. J. Math. 136(2014), 501--550.

%\bibitem{Sta99} R.P. Stanley, Enumerative Combinatorics, vol. 2, in: Cambridge Studies in Advanced Mathematics, No. 62, Cambridge Univ. Press, Cambridge, 1999.

\bibitem{Ste19} J. Stembridge, Sign-twisted Poincar\'e series and odd inversions in Weyl groups, Algebr. Comb. 2(2019), 621--644.

\bibitem{Wac92} M. Wachs, An involution for signed Eulerian numbers, Discrete Math. 99(1992), 59--62.

\end{thebibliography}
\end{document}